# UPPER BOUND FOR THE MOMENT OF SHIFTED VALUES OF CUBIC $L$-FUNCTIONS OVER FUNCTION FIELDS

PRANENDU DARBAR, SAMPA DEY, AND GOPAL MAITI

Abstract. In this paper, we study correlations of shifted values of cubic $L$-functions over function fields and derive an upper bound for moments of these shifted values in the limit where the genus of the corresponding cubic characters tends to infinity over a fixed finite field $\mathbb{F}_q$. Our results apply to the *non-Kummer* case when $q \equiv 2 \pmod 3$. The *Kummer* case, when $q \equiv 1 \pmod 3$, can be treated similarly.

Contents

## 1. Introduction

The goal of this article is to investigate the correlation of shifted values of cubic $L$-functions over the polynomial ring $\mathbb{F}_q[t]$, where $\mathbb{F}_q$ is a finite field. Understanding the correlation of shifted values in families of $L$-functions provides valuable information about their statistical behaviour, including moments and the dependence between values at different shifts along the critical line.

The study of statistical properties of families of $L$-functions is of fundamental importance due to their deep connections with arithmetic objects. In particular, the family of cubic $L$-functions over both number fields and function fields has attracted considerable attention. Various statistical aspects of these families have been investigated in the literature; see, for example, [3, 11, 12, 14, 18, 19, 24, 26].

In comparison with $L$-functions associated with cubic characters, quadratic Dirichlet $L$-functions have been studied much more extensively; see, for example, [1, 10, 15–17, 28, 29]. The study of statistical properties of the cubic family is considerably more challenging due to the more intricate duality and the underlying orthogonality. This is reflected, for instance, in the difficulty of evaluating moments and correlations of shifted values of these $L$-functions, whose conjectural asymptotics are largely guided by predictions from random matrix theory; see [9, 22].

In particular, over function fields, asymptotic formulas are currently known only for the first and second moments of the cubic family [12, 19], whereas the corresponding first four moments have been established for the quadratic family [15–17]. Likewise, significantly stronger non-vanishing results at the central point are available for quadratic Dirichlet $L$-functions in [6] than for the cubic family in [13].

In this article, we will study the moments of the shifted values of cubic $L$-functions over function fields, that helps to understand the correlation among them at different shifts. These objects has been studied in the literature for Riemann zeta function and quadratic Dirichlet $L$-functions over function fields in [8, 10, 25].

### 1.1. Shifted moments and main results.

Before discussing our main results, we first introduce the notation used throughout the paper. Let $\mathbb{F}_q[t]$ denote the polynomial ring over a finite field $\mathbb{F}_q$, where $q$ is an odd prime power. A Dirichlet character modulo $f \in \mathbb{F}_q[t]$ is a multiplicative function on $(\mathbb{F}_q[t]/(f))^*$, extended to $\mathbb{F}_q[t]$ by periodicity and defined to be zero on polynomials that are not coprime to $f$. Let $\mu_3$ denote the group of cubic roots of unity in $\mathbb{C}^*$. A cubic character $\chi$ is a Dirichlet character satisfying $\chi^3 = \chi_0$, where $\chi_0$ denotes the principal character, so that $\chi$ takes values in $\mu_3$. We say that $\chi$ is primitive if its modulus is its conductor, that is, the smallest polynomial modulo which $\chi$ is periodic. We write $L_q(s,\chi)(=\mathcal{L}_q(u,\chi))$ for the associated cubic $L$-function, see Subsection 2.2.3 for more details.

The structure of cubic characters depends on whether the nontrivial cubic roots of unity belong to the base field. If $q \equiv 1 \pmod 3$, the so-called *Kummer* case, then $\mu_3 \subset \mathbb{F}_q$. On the other hand, if $q \equiv 2 \pmod 3$, the *non-Kummer* case, the nontrivial cubic roots of unity are not contained in $\mathbb{F}_q$ but lie in its quadratic extension $\mathbb{F}_{q^2}$.

For the non-Kummer case $q \equiv 2 \pmod 3$, we denote by $\mathcal{C}_g$, the family of primitive cubic Dirichlet characters with conductor of genus $g$. Every character in $\mathcal{C}_g$ is even, since it is trivial on $\mathbb{F}_q^*$. Similarly, when $q \equiv 1 \pmod 3$, we denote by $\widetilde{\mathcal{C}}_g$ the set of *odd* primitive cubic characters of conductor of genus $g$, and these are nontrivial on $\mathbb{F}_q^*$. One could similarly consider the set of all *even* primitive cubic characters; however, the analysis in that case is essentially similar.

Let $e(\theta) = \exp(2\pi i\theta)$. For a fixed integer $m \geq 1$, let $\boldsymbol{k}^{(m)} = (k_1,\ldots,k_m) \in \mathbb{R}^m_{>0}$, and let

$$\boldsymbol{\theta}^{(m)} = (\theta_1,\ldots,\theta_m), \qquad \boldsymbol{\alpha} = (\alpha_1,\ldots,\alpha_m) \in \left[0,\frac{1}{2}\right)^m.$$

We study mixed moments of shifted cubic $L$-functions over the family $\mathcal{C}_g$. More precisely, we define

$$S_g(\boldsymbol{\theta}^{(m)},\boldsymbol{k}^{(m)}) := \sum_{\chi\in\mathcal{C}_g} \left|\mathcal{L}\left(\frac{e(\theta_1)}{q^{1/2+\alpha_1}},\chi\right)\right|^{2k_1} \cdots \left|\mathcal{L}\left(\frac{e(\theta_m)}{q^{1/2+\alpha_m}},\chi\right)\right|^{2k_m}. \tag{1.1}$$

Throughout, we assume that each $\theta_j = \theta_j(g)$ is a real-valued function satisfying for every $i \neq j$, $\lim_{g\to\infty} g|\theta_i-\theta_j|$ also exists or equal to $\infty$. These assumptions naturally describe the microscopic or mesoscopic regime of shifted values.

Our main result establishes an upper bound of the expected order of magnitude for the mixed moments (1.1).

**Theorem 1.1.** *Assume $q \equiv 2 \pmod 3$. Let $\boldsymbol{k}^{(m)}$, $\boldsymbol{\theta}^{(m)}$, and $S_g(\boldsymbol{\theta}^{(m)},\boldsymbol{k}^{(m)})$ be as above, and write $K = \sum_{j=1}^m k_j$. Suppose that $\alpha_j \ll 1/g$ for every $1 \leq j \leq m$. Then, for $g > \exp(4K^2)$,*

$$S_g(\boldsymbol{\theta}^{(m)},\boldsymbol{k}^{(m)}) \ll_{K,q} |\mathcal{C}_g| g^{k_1^2+k_2^2+\cdots+k_m^2} \prod_{1\leq i<j\leq m} \left(\min\left\{\frac{1}{|\theta_i-\theta_j|}, g\right\}\right)^{2k_ik_j}.$$

We give a complete proof in the non-Kummer case, i.e, Theorem 1.1 by combining the methods of and Harper [20] and Soundararajan [30]. The argument for the Kummer case is analogous, and we indicate the necessary modifications whenever they arise (see, for example, Section 2.4). In particular, taking $\boldsymbol{\alpha} = (0,\ldots,0)$, $\boldsymbol{\theta}^{(m)} = (1,\ldots,1)$, Theorem 1.1 recovers the conjectural upper bound predicted in [11, Conjecture 1.2]. The growth rate $g^{k^2}$ agrees with the conjectural order of moments of the Riemann zeta function and is consistent with the unitary symmetry expected for this family. Hence, Theorem 1.1 is comparable to [25, Theorem 1.3]. We also emphasize that the result of Ng–Shen–Wong is conditional on the Riemann Hypothesis, whereas in the function field setting it becomes unconditionally due to Weil [31].

We expect that

$$S_g(\boldsymbol{\theta}^{(m)}, \boldsymbol{k}^{(m)}) \gg_{K,q} |\mathcal{C}_g| g^{k_1+\cdots+k_m^2} \prod_{1\le i<j\le m} \left( \min\left\{ \frac{1}{|\theta_i-\theta_j|}, g \right\} \right)^{2k_ik_j}. \tag{1.2}$$

Establishing the lower bound (1.2) would provide a precise description of the correlation between the values $\mathcal{L}\left(q^{-1/2}e(\theta_1),\chi\right)$ and $\mathcal{L}\left(q^{-1/2}e(\theta_2),\chi\right)$ as the shifts vary. More precisely,

$$\left|\mathcal{L}\left(q^{-1/2}e(\theta_1),\chi\right)\right| \quad \text{and} \quad \left|\mathcal{L}\left(q^{-1/2}e(\theta_2),\chi\right)\right|$$

undergoes a transition when $|\theta_1-\theta_2| \asymp \frac{1}{g}$. In particular, these quantities appear to be asymptotically independent whenever $|\theta_1-\theta_2| \gg \frac{1}{g}$.

Obtaining a matching lower bound (1.2) appears to be substantially more difficult. Even in the case of two shifts, one would require a power-saving error term in the asymptotic formula for the second moment of $\mathcal{L}(q^{-1/2},\chi)$, whereas the strongest result currently available, due to [19], provides only a logarithmic saving.

Finally, let $\mathcal{L}^{(\ell)}(u,\chi)$ denote the $\ell$-th derivative of $\mathcal{L}(u,\chi)$ with respect to $u$. As an immediate consequence of Theorem 1.1, we obtain the following estimate.

**Corollary 1.2.** *Let $\ell, k \in \mathbb{N}$. Then, for $g$ sufficiently large, we have*

$$\sum_{\chi\in\mathcal{C}_g} \left|\mathcal{L}^{(\ell)}\left(q^{-1/2},\chi\right)\right|^{2k} \ll_{k,\ell} |\mathcal{C}_g|\, g^{k^2+2k\ell}.$$

1.2. **Plan of the paper.** The remainder of the paper is organized as follows. In Section 2, we collect the necessary background on the family of cubic characters, the associated $L$-functions, and few lemmas, including orthogonality relations for both the Kummer and non-Kummer cases, which play a crucial role in the proof of our main result.

In Section 3, we establish a Soundararajan-type upper bound, analogous to [30, Corollary A], for the mixed moments $S_g(\boldsymbol{\theta}^{(m)}, \boldsymbol{k}^{(m)})$. This estimate falls short of the conjectured order of magnitude by a factor of $g^{\varepsilon}$, but serves as a key ingredient in the proof of Theorem 1.1. In particular, it enables us to control the contribution from the subfamily of $\mathcal{C}_g$ for which the associated Dirichlet polynomials are large.

In Section 4, we prove Theorem 1.1 by combining the results of the previous sections with Harper's method [20]. The proof proceeds by partitioning the family $\mathcal{C}_g$ into several subfamilies according to the size of the associated Dirichlet polynomials and estimating the contribution from each subfamily separately. Finally in Section 5, we prove Corollary 1.2.

## 2. Preliminary Lemmas

We begin this section with some preliminaries on function fields, basic properties of cubic $L$-functions and their spectral interpretation. Most of the definitions and results are from [27].

2.1. **Basic facts on $\mathbb{F}_q[t]$.** We start by fixing a finite field $\mathbb{F}_q$ of cardinality $q = p^r$, $r \geq 1$ for a $p$ such that $\mathrm{char}(\mathbb{F}_q) \neq 2, 3$. We denote the polynomial ring over $\mathbb{F}_q$ as $\mathbb{A} := \mathbb{F}_q[t]$. For a polynomial $f$ in $\mathbb{F}_q[t]$, it's degree will be denoted by either $\deg(f)$ or $d(f)$.

The set of all monic square-free, monic polynomials and monic irreducible polynomials of degree $n$ are denoted by $\mathcal{H}_{q,n}, \mathcal{M}_{n,q}$ (or simply $\mathcal{M}_n, \mathcal{H}_n$ as we fix $q$) and $\mathcal{P}_{n,q}$ (or simply $\mathcal{P}_n$) respectively. Let $\mathcal{M} = \cup_{n\geq 1}\mathcal{M}_n$ and $\mathcal{P} = \cup_{n\geq 1}\mathcal{P}_n$. We also denote the set of all monic polynomials and monic irreducible polynomials of degree less or equal to $n$ by $\mathcal{M}_{\leq n,q}$ (or simply $\mathcal{M}_{\leq n}$) and $\mathcal{P}_{\leq n,q}$ (or simply $\mathcal{P}_{\leq n}$) respectively. Let $\mathcal{H}_{q,\leq n}$ be the set of monic square-free polynomials of degree less than or equal to $n$, and set $\mathcal{H}_q = \cup_{n\geq 1}\mathcal{H}_n$.

Observe that for $n \geq 1$, $|\mathcal{M}_n| = q^n$. For any $f \in \mathbb{F}_q[t] \setminus \{0\}$, the norm of $f$ is defined by $|f| = q^{d(f)}$. If $f = 0$, set $|f| = 0$. The function field analogue of the prime number theorem [27, Theorem 2.2] asserts that

$$|\mathcal{P}_{n,q}| = \frac{q^n}{n} + O\Big(\frac{q^{\frac{n}{2}}}{n}\Big). \tag{2.1}$$

The zeta function of $\mathbb{A}$, denoted by $\zeta_{\mathbb{A}}(s)$ is defined by

$$\zeta_q(s) = \sum_{f\in\mathcal{M}} \frac{1}{|f|^s} = \prod_{P\in\mathcal{P}} \left(1 - |P|^{-s}\right)^{-1}, \qquad \mathrm{Re}(s) > 1.$$

It can be shown that $\zeta_q(s) = \frac{1}{1-q^{1-s}}$, and this provides an analytic continuation of zeta function to the complex plane with a simple pole at $s = 1+(2\pi i n)/\ln q$, where $n \in \mathbb{Z}$. Using the change of variable $u = q^{-s}$, we define

$$\mathcal{Z}_q(u) = \sum_{f\in\mathcal{M}} u^{d(f)} = \frac{1}{1-qu}, \quad |u| < \frac{1}{q}. \tag{2.2}$$

2.2. **Cubic Dirichlet characters and properties of their $L$-functions.**

2.2.1. *Kummer setting.* First we consider the Kummer case, that is $q \equiv 1 \pmod 3$. So the cubic roots of 1 lie in $\mathbb{F}_q^*$. We fix an isomorphism $\Omega$ between the sets of cubic roots of unity $\mu_3 \subset \mathbb{C}^*$ and the cubic roots of unity in $\mathbb{F}_q^*$. For a prime polynomial $P \in \mathbb{F}_q[t]$ and $f \in \mathbb{F}_q[t]$ with $P \nmid f$, there exists a unique $\alpha \in \mu_3$ such that

$$f^{\frac{|P|-1}{3}} \equiv \Omega(\alpha) \pmod P.$$

We set the cubic character $\chi_P$ modulo $P$ by

$$\chi_P(f) = \alpha.$$

Note that, depending on the choices of the isomorphism $\Omega$, there are two such cubic characters $\chi_P$ and $\overline{\chi_P} = \chi_P^2$. This construction can be lift to any monic polynomial in $\mathbb{F}_q[t]$ by multiplicativity. The cubic character is primitive modulo $F$ if for any $P^e|F$ either $e = 1$ or 2. For any character $\chi$ on $\mathbb{F}_q[t]$, if it is trivial on $\mathbb{F}_q^*$ then we call it *even* character and *odd* otherwise. It follows that the conductor of the primitive cubic characters are the square-free $F \in \mathcal{M}$, and for each such conductor there are $2^{\omega(F)}$ characters, where $\omega(F)$ is the number of irreducible dividing $F$.

2.2.2. *Non-Kummer setting.* When $q \equiv 2 \pmod 3$, the cubic characters are basically a subset of the cubic characters over $\mathbb{F}_{q^2}[t]$. Precisely, in this case, the cubic character $\chi_P$ can be constructed as in the Kummer setting if and only if $P$ is of even degree. The conductor of the primitive cubic characters are the square-free polynomials $F \in \mathbb{F}_q[t]$ supported on irreducibles of even degree. Also for each such conductor, there are $2^{\omega(F)}$ characters. We found detailed discussion in [4, 11, 12].

Moreover, if $q \equiv 1 \pmod 6$, then there is perfect cubic reciprocity which says, for $f, g \in \mathbb{F}_q[t]$, monic such that $(f, g) = 1$, the corresponding cubic residue symbol $\chi_a$ and $\chi_b$ satisfies

$$\chi_a(b) = \chi_b(a).$$

For the general cubic reciprocity one can see [27, Theorems 3.3 and 3.5].

2.2.3. *Cubic $L$-functions.* The $L$-function associated to a primitive cubic character $\chi$ of conductor $F$ is defined by

$$L_q(s, \chi_F) = \sum_{f \in \mathcal{M}} \frac{\chi_F(f)}{|f|^s} = \prod_{\substack{P \in \mathcal{P} \\ P \nmid F}} \left(1 - \chi_F(P)\,|P|^{-s}\right)^{-1}, \quad \mathrm{Re}(s) > 1.$$

Using the change of variable $u = q^{-s}$, we have

$$\mathcal{L}_q(u, \chi_F) = \sum_{f \in \mathcal{M}} \chi_F(f)\, u^{d(f)} = \prod_{\substack{P \in \mathcal{P} \\ P \nmid F}} \left(1 - \chi_F(P)\, u^{d(P)}\right)^{-1}, \quad |u| < \frac{1}{q}.$$

By the orthogonality relations of $\chi$ [cf. [27], Proposition 4.3], we see that if $n \geq \deg(F)$ then

$$\sum_{f \in \mathcal{M}_n} \chi_F(f) = 0.$$

It implies that $\mathcal{L}_q(u, \chi_F)$ is a polynomial of degree at most $\deg(F) - 1$.

2.3. **Spectral Interpretation.** Let $C$ be a non-singular projective curve of genus $g$ over $\mathbb{F}_q(t)$ whose function field is a cyclic extension of $\mathbb{F}_q(t)$. For each extension field $\mathbb{F}_{q^k}$ on $\mathbb{F}_q$ of degree $k$, denote by $N_k(C)$ the number of points of $C$ in $\mathbb{F}_{q^k}$. The zeta function associated to $C$ defined as

$$Z_C(u) = \exp\left(\sum_{k=1}^{\infty} N_k(C) \frac{u^k}{k}\right), \quad |u| < \frac{1}{q}.$$

From the Weil conjectures, it is known to be a rational function of $u$ of the form

$$Z_C(u) = \frac{P_C(u)}{(1-u)(1-qu)}.$$

Additionally, we know that $P_C(u)$ is a polynomial of degree $2g$ with integer coefficients, satisfying a functional equation

$$P_C(u) = (qu^2)^g P_C\left(\frac{1}{qu}\right).$$

The Riemann Hypothesis, proved by Weil [31], says that the zeros of $P_C(u)$ all lie on the circle $|u| = \frac{1}{\sqrt{q}}$. Thus one may give a spectral interpretation of $P_C(u)$ as the characteristic polynomial of a $2g \times 2g$ unitary matrix $\Theta_C$:

$$P_C(u) = \det\left(I - u\sqrt{q}\Theta_C\right).$$

The eigenvalues $e^{2\pi i \theta_j}$ of $\Theta_C$ correspond to the zeros, $q^{-1/2} e^{-2\pi i \theta_j}$, of $Z_C(u)$. The matrix $\Theta_C$ is called the unitarized Frobenius class of $C$.

To put this in the context of our case,

$$P_C(u) = \mathcal{L}(u, \chi_F)\mathcal{L}(u, \overline{\chi_F}),$$

where $\chi_F$ and $\overline{\chi_F}$ are two cubic Dirichlet characters of conductor $F$ of the function fields of $C$. Considering the prime at infinity over $\mathbb{F}_q(t)$, we have

$$\mathcal{L}_C(u, \chi) = \begin{cases} \mathcal{L}_q(u, \chi) & \text{if } \chi \text{ is odd,} \\ \frac{\mathcal{L}_q(u,\chi)}{1-u} & \text{if } \chi \text{ is even.} \end{cases}$$

From the Riemann–Hurwitz formula,

$$\deg(F) = g + 2 - \begin{cases} 0 & \text{if } \chi \text{ is even,} \\ 1 & \text{if } \chi \text{ is odd.} \end{cases}$$

2.4. **Orthogonality of characters over the families.** The family of primitive cubic Dirichlet characters with conductor of genus $g$ depending on Kummer or non-Kummer setting. Suppose $q$ is odd and $q \equiv 1 \pmod 3$. Define

$$\widetilde{\mathcal{C}}_g = \{\chi : \chi \text{ is an } \textit{odd} \text{ primitive character with conductor of genus } g\} \tag{2.3}$$

In this case, from [12, Lemma 2.10], for every $\varepsilon > 0$, we have

$$\#\widetilde{\mathcal{C}}_g = B_1 g q^g + B_2 q^g + O\left(q^{\frac{g}{2}(1+\varepsilon)}\right), \tag{2.4}$$

where $B_1 = q\mathcal{G}(1/q), B_2 = (2q\mathcal{G}(1/q) - \mathcal{G}'(1/q))$, and

$$\mathcal{G}(u) = \prod_{P \in \mathcal{P}} \left(1 - 3u^{2d(P)} + 2u^{3d(P)}\right).$$

Suppose that $q$ is odd, and $q \equiv 2 \pmod 3$. Define

$$\mathcal{C}_g = \{\chi : \chi \text{ is an } \textit{even} \text{ primitive character with conductor of genus } g\}. \tag{2.5}$$

In this case, we must have $g$ is *even*, otherwise the set becomes empty. From [12, Lemmma 2.10], we have

$$\#\mathcal{C}_g = B_3 q^{g+2} + O(q^{\frac{g}{2}(1+\epsilon)}), \tag{2.6}$$

where

$$B_3 = \prod_{\substack{P \in \mathcal{P} \\ d(P) \text{ odd}}} \left(1 - u^{2d(P)}\right) \prod_{\substack{P \in \mathcal{P} \\ d(P) \text{ even}}} \left(1 - 3u^{2d(P)} + 2u^{3d(P)}\right).$$

**Lemma 2.1** (Character sum over a cube)**.** *Let $q \equiv 2 \pmod 3$. For $f \in \mathcal{M}$,*

$$\sum_{\chi \in \mathcal{C}_g} \chi(f^3) = (\#\mathcal{C}_g) \prod_{\substack{2|\deg(P) \\ P|f}} \left(\frac{|P|}{|P|+2}\right) + O\left(q^{(1/2+\epsilon)g} 2^{\omega(f)}\right)$$

*where $\omega(f)$ counts the number of distinct prime factors of $f$.*

*Let $q \equiv 1 \pmod 3$. Then we have*

$$\sum_{\chi \in \widetilde{\mathcal{C}}_g} \chi(f^3) = (\#\widetilde{\mathcal{C}}_g) \prod_{P|f} \left(\frac{|P|}{|P|+2}\right) + O\left(q^{(1/2+\epsilon)g} 2^{\omega(f)}\right).$$

*Proof.* We consider the non-Kummer case. The proof for the Kummer case follows along the same lines. We can write

$$\sum_{\chi \in \mathcal{C}_g} \chi(f^3) = \sum_{\substack{\chi \in \mathcal{C}_g \\ (\mathrm{Cond}(\chi), f) = 1}} 1 = \sum_{\substack{F \in \mathcal{M} \\ (F, f) = 1}} a_g(F),$$

where $a_g(F)$ counts the number of cubic primitive characters of conductor $F$ and degree $(g+2)$. By Perron's formula over function fields [15, Equation (5.3)], we have

$$\sum_{\chi \in \mathcal{C}_g} \chi(f^3) = \frac{1}{2\pi i} \int_{|u|=q^{-2}} \mathcal{A}(u) \frac{du}{u^{g+3}},$$

where

$$\mathcal{A}(u) = \sum_{\substack{F \in \mathcal{M} \\ (F,f)=1}} a(F) u^{d(F)} = \prod_{\substack{2|d(P) \\ P|f}} \left(1 + 2u^{d(P)}\right)^{-1} \prod_{2|d(P)} \left(1 + 2u^{d(P)}\right).$$

Notice that this power series is analytic for $|u| < q^{-1}$ with simple poles at $u = \pm q^{-1}$. We know that the zeta function $\mathcal{Z}_{q^2}(u)$ has simple poles at those points. So we can multiply $\mathcal{A}(u)$ by $\mathcal{Z}_{q^2}(u)$ to get an series which is analytic in a wider circular region. This means by defining

$$\begin{aligned}
\mathcal{G}(f;u) &= \mathcal{A}(u)(1-qu)(1+qu) \\
&= \prod_{\substack{2|d(P) \\ P|f}} \left(1 + 2u^{d(P)}\right)^{-1} \prod_{2|d(P)} \left(1 + 2u^{d(P)}\right) \prod_{P} \left(1 - u^{d(P)}\right) \left(1 - (-1)^{d(P)} u^{d(P)}\right) \\
&= \mathcal{G}(u) \prod_{\substack{2|d(P) \\ P|f}} \left(\frac{|P|}{|P|+2}\right),
\end{aligned}$$

where $\mathcal{G}(u)$ is analytic for $|u| < q^{-1/2}$. We move the contour of integration from $|u| = q^{-2}$ to $|u| = q^{-1/2-\varepsilon}$ by encountering simple poles at $u = \pm q^{-1}$, we obtain

$$\begin{aligned}
&\sum_{\chi \in \mathcal{C}_g} \chi(f^3) = \frac{1}{2\pi i} \int_{|u|=q^{-2}} \frac{\mathcal{G}(f;u)}{(1-qu)(1+qu)} \frac{du}{u^{g+3}} \\
&= \left(\frac{\mathcal{G}(1/q)}{2} + (-1)^g \frac{\mathcal{G}(1/q)}{2}\right) q^{g+2} \prod_{\substack{2|d(P) \\ P|f}} \left(\frac{|P|}{|P|+2}\right) + O\left(2^{\omega(f)} q^{(1/2+\varepsilon)g}\right).
\end{aligned}$$

Since $g$ is even, the result follows. □

**Lemma 2.2** (Pólya–Vinogradov inequality for non-Kummer extension)**.** *Let $q \equiv 2 \pmod 3$. For $f \in \mathcal{M}$ not a perfect cube, write $f = f_1 f_2^2 f_3^3$, where $f_1, f_2$ are square-free and co-prime. Then for any $\varepsilon > 0$,*

$$\sum_{\chi \in \mathcal{C}_g} \chi(f) \ \ll_\varepsilon \ g q^{g/2} |f_1 f_2^2|^\varepsilon \prod_{P|f_3} \left(1 + \frac{1}{|P|^{1/2}}\right).$$

*Proof.* Using the definition of $\mathcal{C}_g$, we can express

$$\sum_{\chi\in\mathcal{C}_g}\chi(f)=\sum_{\substack{F\in\mathcal{H}_{q^2,g/2+1}\\ P|F\implies P\notin\mathbb{F}_q[t]}}\chi_F(f)=\sum_{\substack{h\in\mathcal{M}_{q,\leq g/2+1}\\ (h,f)=1}}\mu(h)\sum_{\substack{h\in\mathcal{H}_{q^2,\leq g/2+1-d(h)}\\ (F,h)=1}}\chi_F(f),$$

where we use the detector

$$\sum_{\substack{h\in\mathbb{F}_q[t]\\ h|F}}\mu(h)=\begin{cases}1 & F \text{ has no prime divisor in } \mathbb{F}_q[t],\\ 0 & \text{otherwise.}\end{cases}$$

Using Perron's formula, we estimate the inner sum as

$$\sum_{\substack{h\in\mathcal{H}_{q^2,\leq g/2+1-d(h)}\\ (F,h)=1}}\chi_F(f)=\frac{1}{2\pi i}\int_{|u|=q^{-3}}\mathcal{B}(f;u)\frac{du}{u^{g/2+2-d(h)}},$$

where we use cubic reciprocity to simplify

$$\begin{aligned}\mathcal{B}(f;u)&=\sum_{\substack{F\in\mathcal{H}_{q^2}\\ (F,h)=1}}\chi_f(F)u^{\deg(f)}=\prod_{\substack{P\in\mathbb{F}_{q^2}[t]\\ P\nmid fh}}\left(1+\chi_f(P)u^{\deg(P)}\right)\\ &=\frac{\mathcal{L}_{q^2}(u,\chi_f)}{\mathcal{L}_{q^2}(u^2,\overline{\chi}_f)}\prod_{\substack{P\in\mathbb{F}_{q^2}[t]\\ P\nmid f;P|h}}\frac{1-\chi_f(P)u^{d(P)}}{1-\overline{\chi}_f(P)u^{2d(P)}}.\end{aligned}$$

On the other hand, we can write

$$\mathcal{L}_{q^2}(u,\chi_f)=\mathcal{L}_{q^2}(u,\chi_{f_1f_2^2})\prod_{\substack{P\in\mathbb{F}_{q^2}[t]\\ P\nmid f_1f_2\\ P|f_3}}\left(1-\chi_{f_1f_2^2}(P)u^{d(P)}\right).$$

We use [5, Theorem 5.1] to obtain that, for $|u|\leq q^{-1}$,

$$\mathcal{L}_{q^2}(u,\chi_f)\ll_\varepsilon|f_1f_2^2|^{2\varepsilon}\prod_{P|f_3}\left(1+\frac{1}{|P|^{1/2}}\right).$$

Furthermore, using [12, Lemma 2.7], for $|u|\leq q^{-1}$, we bound

$$\mathcal{L}_{q^2}(u^2,\overline{\chi}_f)\gg_\varepsilon|f_1f_2^2|^{-2\varepsilon}.$$

By moving the contour of integration to $|u|=q^{-1}$, from the above bounds, we obtain

$$\sum_{\substack{h\in\mathcal{H}_{q^2,\leq g/2+1-d(h)}\\ (F,h)=1}}\chi_F(f)\ll_\varepsilon q^{g/2-\deg(h)}|f_1f_2^2|^{4\varepsilon}\prod_{P|f_3}\left(1+\frac{1}{|P|^{1/2}}\right)\prod_{P|h}\left(1+\frac{1}{|P|^{1/2}}\right).$$

Therefore, finally we conclude that

$$\sum_{\chi\in\mathcal{C}_g}\chi(f) \ll q^{g/2}|f_1f_2^2|^{4\varepsilon}\prod_{P|f_3}\left(1+\frac{1}{|P|^{1/2}}\right)\sum_{h\in\mathcal{M}_{q,\le g/2+1}}\frac{1}{|h|}\prod_{P|h}\left(1+\frac{1}{|P|^{1/2}}\right)$$
$$\ll q^{g/2}g|f_1f_2^2|^{4\varepsilon}\prod_{P|f_3}\left(1+\frac{1}{|P|^{1/2}}\right),$$

which finishes the proof. □

**Lemma 2.3** (Pólya–Vinogradov inequality for Kummer extension)**.** *Let $q\equiv 1 \pmod 3$. Then for any $\varepsilon>0$, when $f\in\mathcal{M}$ not a perfect cube, we have*

$$\sum_{\chi\in\widetilde{\mathcal{C}}_g}\chi(f) \ll_\varepsilon q^{g/2}g^7|f|^\varepsilon.$$

*Proof.* Using the definition of $\widetilde{\mathcal{C}}_g$, for each primitive cubic character $\chi_{F_1F_2^2}$, we have that for $\alpha\in\mathbb{F}_q^*$,

$$\chi_{F_1F_2^2}(\alpha)=\Omega^{-1}\left(\alpha^{\frac{q-1}{3}(\deg(F_1)+2\deg(F_2))}\right).$$

Moreover, since $\chi_{F_1F_2^2}$ is odd, the restriction to $\mathbb{F}_q^*$ is $\chi_3$ when $\deg(F_1)+2\deg(F_2)\equiv 1 \pmod 3$, and $\chi_3^2$ when $\deg(F_1)+2\deg(F_2)\equiv 2 \pmod 3$, where $\chi_3$ is a cubic character over $\mathbb{F}_q^*$. Then, since the conductor of $\chi_{F_1F_2^2}$ is $F=F_1F_2$,

$$\deg(F_1)+\deg(F_2)=g+1.$$

Therefore, we can express

$$\sum_{\chi\in\widetilde{\mathcal{C}}_g}\chi(f)=\sum_{\substack{d_1+d_2=g+1\\ d_1+2d_2\equiv 1 \pmod 3}}\sum_{\substack{F_1\in\mathcal{H}_{q,d_1}\\ F_2\in\mathcal{H}_{q,d_2}\\ (F_1,F_2)=1}}\chi_{F_1F_2^2}(f),$$

Using cubic reciprocity and [12, Corollary 5.5], we have

$$\sum_{\substack{F_1\in\mathcal{H}_{q,d_1}\\ F_2\in\mathcal{H}_{q,d_2}\\ (F_1,F_2)=1}}\chi_{F_1F_2^2}(f)=\sum_{\substack{H\in\mathcal{M}_{q,\le\min\{d_1,d_2\}}\\ (H,f)=1}}\mu(H)\sum_{\substack{R_1\in\mathcal{M}_{q,\le d_1-d(H)}\\ R_1|H}}\mu(R_1)\chi_f(R_1)$$
$$\times\sum_{\substack{R_2\in\mathcal{M}_{q,\le d_2-d(H)}\\ R_2|H}}\mu(R_2)\chi_f(R_2)^2\sum_{\substack{D_1\in\mathcal{M}_{q,\le\frac{d_1-d(H)-d(R_1)}{2}}\\ (D_1,R_1)=1}}\mu(D_1)\chi_f(D_1)^2$$
$$\times\sum_{\substack{D_2\in\mathcal{M}_{q,\le\frac{d_2-d(H)-d(R_2)}{2}}\\ (D_2,R_2)=1}}\mu(D_2)\chi_f(D_2)\sum_{\substack{L_1\in\mathcal{M}_{q,d_1-d(D_1^2HR_1)}\\ L_2\in\mathcal{M}_{q,d_2-d(D_2^2HR_2)}}}\chi_f(L_1)\overline{\chi_f}(L_2).$$

Next using Perron's formula for the sums over $L_j, j = 1, 2$, we obtain

$$\sum_{d(L_j)=d_j-d(D_j^2HR_j)} \chi_f(L_j) = \frac{1}{2\pi i}\int_{|u|=q^{-1/2}} \frac{\mathcal{L}_q(u,\chi_f)}{u^{d_j-d(D_j^2HR_j)}}\frac{du}{u}.$$

Since $f$ is not a cube, the cubic $L$-functions associated to non trivial character $\chi_f$ in the above expression have no poles and we can shift contour of integration on the circle $|u| = q^{-1/2}$. Using the Lindelöf hypothesis [12, Lemma 2.6]), on $|u| = q^{-1/2}$, we obtain

$$|\mathcal{L}_q(u,\chi_f)| \ll_\varepsilon |f|^\varepsilon.$$

Therefore, we have

$$\sum_{d(L_j)=d_j-d(D_j^2HR_j)} \chi_f(L_j) \ll |f|^\varepsilon \frac{q^{d_j/2}}{|D_j^2HR_j|^{1/2}}.$$

Putting together all these bound, we have

$$\sum_{\chi\in\widetilde{\mathcal{C}}_g} \chi(f) \ll q^{\frac{g+1}{2}}|f|^\varepsilon \sum_{\substack{d_1+d_2=g+1\\ d_1+2d_2\equiv 1 \ (\mathrm{mod}\ 3)}} \sum_{\substack{H\in\mathcal{M}_{q,\le\min\{d_1,d_2\}}\\ (H,f)=1}} \frac{1}{|H|} \sum_{\substack{R_1\in\mathcal{M}_{q,\le d_1-d(H)}\\ R_1|H}} \frac{1}{|R_1|^{1/2}}$$

$$\times \sum_{\substack{R_2\in\mathcal{M}_{q,\le d_2-d(H)}\\ R_2|H}} \frac{1}{|R_2|^{1/2}} \sum_{\substack{D_1\in\mathcal{M}_{q,\le\frac{d_1-d(H)-d(R_1)}{2}}\\ (D_1,R_1)=1}} \frac{1}{|D_1|} \sum_{\substack{D_2\in\mathcal{M}_{q,\le\frac{d_2-d(H)-d(R_2)}{2}}\\ (D_2,R_2)=1}} \frac{1}{|D_2|}.$$

Note that, the sums over $D_1$ and $D_2$ are bounded by $d_1$ and $d_2$ respectively. Also, the sums over $R_1$ and $R_2$ are both bounded by $\tau(H)$, where $\tau(H)$ is the divisor function of $H$. Therefore, we obtain

$$\sum_{\chi\in\widetilde{\mathcal{C}}_g} \chi(f) \ll q^{\frac{g}{2}}|f|^\varepsilon \sum_{\substack{d_1+d_2=g+1\\ d_1+2d_2\equiv 1 \ (\mathrm{mod}\ 3)}} d_1d_2(\min\{d_1,d_2\})^4 \ll g^7q^{\frac{g}{2}}|f|^\varepsilon.$$

This finishes the proof. □

2.5. **Bounding logarithm of $L$-functions by a Dirichlet polynomial.** The following lemma bounds logarithm of the cubic $L$-functions by a short Dirichlet polynomial and the proof follows from [7, Proposition 4.3]. From this section onwards, ln and $\log_q$ denote the natural logarithm and the logarithm to base $q$, respectively.

**Lemma 2.4.** *Let $\chi$ be a cubic Dirichlet character with conductor of genus $g$ over $\mathbb{F}_q[t]$. Assume that $\alpha \ll 1/g$ and $\theta \in [0,\pi)$. Then for $N \le g+2$, we have*

$$\ln\left|\mathcal{L}\left(q^{-\frac{1}{2}-\alpha}e(\theta),\chi\right)\right| \le \mathrm{Re}\left(\sum_{\deg(f)\le N} \frac{\Lambda(f)\chi(f)(N-\deg(f))}{N\deg(f)|f|^{\frac{1}{2}+\alpha-\frac{2\pi i\theta}{\ln q}+\frac{1}{N\ln q}}}\right)+\frac{g+2}{N}+O(1).$$

## 3. A Soundararajan type upper bound

Consider $q \equiv 2 \pmod 3$. In this section, we prove the following proposition by adapting the method of [30, Corollary A]. This result serves as one of the key ingredient in establishing the sharp upper bound of Theorem 1.1.

**Proposition 3.1.** *With the assumption as in the Theorem 1.1, for any $\epsilon > 0$,*

$$S_g(\boldsymbol{\theta}^{(m)}, \boldsymbol{k}^{(m)}) \ll_{\boldsymbol{k}^{(m)},\varepsilon} |\mathcal{C}_g| g^{\varepsilon} \exp\left(\sigma\left(\boldsymbol{\theta}^{(m)}, g\right)\right),$$

*where*

$$\sigma(\boldsymbol{\theta}^{(m)}, g) = \left(\sum_{j=1}^{m} k_j^2\right) \ln g \, + \, 2\sum_{i<j} k_i k_j \left(\ln\left(\min\left\{\tfrac{1}{|\theta_i-\theta_j|}, g\right\}\right)\right). \tag{3.1}$$

We need the following lemma, which provides an upper bound for the moments of a certain Dirichlet polynomial over the family.

**Lemma 3.2.** *Let $\ell$ and $y$ be integers such that $6\ell y \le g+2$. For any complex numbers $a(P)$ with $|a(P)| \ll 1$, we have*

$$\sum_{\chi\in\mathcal{C}_g} \left| \sum_{\deg(P)\le y} \frac{\chi(P)a(P)}{\sqrt{P}} \right|^{2\ell} \ll q^g \frac{(\ell!)^2 5^{2\ell/3}}{\lfloor 2\ell/3\rfloor! 9^{\ell/3}} \left( \sum_{\deg(P)\le y} \frac{|a(P)|^2}{|P|} \right)^{\ell}.$$

*Proof.* See [11, Lemma 6.2]. □

### 3.1. Proof of Proposition 3.1.

To keep things simple we use the notation $\boldsymbol{\theta}$ instead of $\boldsymbol{\theta}^{(m)}$. Consider

$$\Psi_g(\boldsymbol{\theta}, V) = \#\left\{ \chi \in \mathcal{C}_g : \sum_{j=1}^{m} 2k_j \ln\left|\mathcal{L}\Big(\frac{e(\theta_j)}{q^{\frac{1}{2}+\alpha_j}}, \chi\Big)\right| \ge V \right\},$$

for sufficiently large $g$ and for all $V > 2$. Then we can express

$$S_g(\boldsymbol{\theta}^{(m)}, \boldsymbol{k}^{(m)}) = \int_{-\infty}^{\infty} \Psi_g(\boldsymbol{\theta}, V) \exp\big(V\big) dV. \tag{3.2}$$

To bound $S_g(\boldsymbol{\theta}^{(m)}, \boldsymbol{k}^{(m)})$, we will estimate an upper bound of $\Psi_g(\boldsymbol{\theta}, V)$ for different ranges of $V$. The Lemma 2.4 leads us

$$
\begin{aligned}
2\sum_{j=1}^{m} k_j \ln\left|\mathcal{L}\Big(\frac{e(\theta_j)}{q^{\frac{1}{2}+\alpha_j}}, \chi\Big)\right| \leq\ & 2\,\mathrm{Re}\left(\sum_{\deg(P)\leq N}\sum_{j=1}^{m} k_j \frac{\chi(P)(N-\deg(P))}{N|P|^{\frac{1}{2}+\alpha_j-\frac{2\pi i\theta_j}{\ln q}+\frac{1}{N\ln q}}}\right)\\
&+2\,\mathrm{Re}\left(\sum_{\deg(P)\leq N/2}\sum_{j=1}^{m} k_j \frac{\chi(P^2)(N-2\deg(P))}{2N|P|^{1+2\alpha_j-\frac{4\pi i\theta_j}{\ln q}+\frac{2}{N\ln q}}}\right)\\
&+2\sum_{\ell\geq 3}\mathrm{Re}\left(\sum_{\deg(P)\leq N/\ell}\sum_{j=1}^{m} k_j \frac{\chi(P^\ell)(N-\ell\deg(P))}{\ell N|P|^{\ell\left(\frac{1}{2}+\alpha_j-\frac{2\pi i\theta_j}{\ln q}+\frac{1}{N\ln q}\right)}}\right)\\
&+\frac{2K(g+2)}{N}+O_K(1),
\end{aligned}
\tag{3.3}
$$

where we recall $K=\sum_{j=1}^{m} k_j$. For a non-trivial cubic character $\chi\neq\chi_0$ with conductor of genus $g$, we have ([cf. [21, Eq. (3.33)])

$$
\sum_{\deg P=n}\chi(P^2)=O\left(\frac{q^{n/2}g}{n}\right). \tag{3.4}
$$

Note that

$$
\begin{aligned}
&\sum_{\deg(P)\leq N/2}\sum_{j=1}^{m} k_j \frac{\chi(P^2)(N-2\deg(P))}{2N|P|^{1+2\alpha_j-\frac{4\pi i\theta_j}{\ln q}+\frac{2}{N\ln q}}}\\
&\quad=\sum_{j=1}^{m}\frac{k_j}{2N}\sum_{n=1}^{N/2}\frac{N-2n}{q^{n\left(1+2\alpha_j+\frac{2}{N\ln q}-\frac{4\pi i\theta_j}{\ln q}\right)}}\sum_{\deg(P)=n}\chi(P^2)\\
&\quad=\sum_{j=1}^{m}\frac{k_j}{2N}\left(\sum_{n=1}^{\lfloor 3\log_q g\rfloor}+\sum_{\lfloor 3\log_q g\rfloor<n}^{N/2}\right)\frac{N-2n}{q^{n\left(1+2\alpha_j+\frac{2}{N\ln q}-\frac{4\pi i\theta_j}{\ln q}\right)}}\sum_{\deg(P)=n}\chi(P^2).
\end{aligned}
$$

Then the first sum is trivially bounded by $O_K(\ln\log_q g)$ by using the prime number theorem (2.1). By partial summation and using equation (3.4), the second term is $O_K(1)$. Therefore, the contributions from the second and third terms to the right-hand side of the inequality (3.3) are $O_K(\ln\log_q g)$ and $O_K(1)$, respectively.

Let us consider

$$h_1(P) = \mathrm{Re}\sum_{j=1}^{m} k_j |P|^{-\alpha_j + \frac{2\pi i\theta_j}{\ln q} + \frac{1}{N\ln q}}, \tag{3.5}$$

and $a(P) = 1 - \frac{\deg P}{N}$. Certainly $a(P) \ll 1$, since $\deg(P) \geq 1$.

Therefore, we have

$$\sum_{j=1}^{m} 2k_j \ln\left|\mathcal{L}\Big(\frac{e(\theta_j)}{q^{\frac{1}{2}+\alpha_j}}, \chi\Big)\right| \leq \mathcal{T}_1(\chi) + \mathcal{T}_2(\chi) + \frac{(2g+4)K}{N} + O_K(\ln\log_q g),$$

where

$$\mathcal{T}_1(\chi) = 2\sum_{d(P)\leq N_0} \frac{\chi(P)h_1(P)a(P)}{|P|^{1/2}}$$

and

$$\mathcal{T}_2(\chi) = 2\sum_{N_0 < d(P) \leq N} \frac{\chi(P)h_1(P)a(P)}{|P|^{1/2}}.$$

From now onward, for the sake of simplicity, we write $\sigma(\boldsymbol{\theta}, g)$ simply as $\boldsymbol{\sigma}$. We consider the following ranges of $V$. For the range $-\infty < V \leq \sqrt{\ln g}$, the contribution to the integral is

$$\int_{-\infty}^{\infty} \Psi_g(\boldsymbol{\theta}, V)\exp\big(V\big)dV \ll |\mathcal{C}_g|\exp\left(\sqrt{\ln g}\right) \ll |\mathcal{C}_g| g^{o(1)}.$$

Therefore, it suffices to prove Proposition 3.1, with the assumption that $\sqrt{\ln g} \leq V \leq \frac{Kg}{\log_q g}$. Let us define the quantity $A$ as

$$A = \begin{cases} \frac{\ln\boldsymbol{\sigma}}{2}, & \text{if } \sqrt{\ln g} \leq V \leq \boldsymbol{\sigma}, \\ \frac{\boldsymbol{\sigma}\ln\boldsymbol{\sigma}}{2V}, & \text{if } \boldsymbol{\sigma} < V \leq \frac{\boldsymbol{\sigma}\ln\boldsymbol{\sigma}}{25K}, \\ 7K, & \text{if } V > \frac{\boldsymbol{\sigma}\ln\boldsymbol{\sigma}}{25K}. \end{cases}$$

Take $N = \frac{(2g+4)A}{V}$ and $N_0 = \frac{N}{\log_q g}$,. If $\chi \in \Psi_g(\boldsymbol{\theta}, V)$ then we must have either

$$\mathcal{T}_1(\chi) \geq V\left(1 - \frac{6K}{A}\right) := V_1$$

for $g$ large enough since $\sqrt{\ln g} < V$, or

$$\mathcal{T}_2(\chi) \geq \frac{KV}{A} := V_2.$$

For $i = 1, 2$, consider

$$\Psi_g(\boldsymbol{\theta}, V_i) = \#\left\{\chi \in \mathcal{C}_g : \mathcal{T}_i(\chi) \geq V_i\right\},$$

Using Lemma 3.2, we obtain

$$\sum_{\chi\in\mathcal{C}_g} |\mathcal{T}_2(\chi)|^{2\ell} \ll |\mathcal{C}_g| \frac{(\ell!)^2 5^{2\ell/3}}{\lfloor 2\ell/3\rfloor ! 9^{\ell/3}} \left( \sum_{N_0<d(P)\leq N} \frac{|h_1(P)|^2 |a(P)|^2}{|P|} \right)^{\ell}$$

for any $\ell$ such that $6\ell N \leq g+2$, which implies that $\ell \leq \frac{g+2}{6N} \leq \frac{V}{12A}$. Therefore, by taking $\ell = \lfloor \frac{V}{12A} \rfloor$ and using Markov's inequality and Stirling's formula, it follows that

$$\begin{aligned}
\Psi_g(\boldsymbol{\theta}, V_2) &\leq \frac{1}{{V_2}^{2\ell}} \Big( \sum_{\chi\in\mathcal{C}_g} |\mathcal{T}_2(\chi)|^{2\ell} \Big) \\
&\ll |\mathcal{C}_g| \Big( \frac{A}{KV} \Big)^{2\ell} \frac{(\ell!)^2 5^{2\ell/3}}{\lfloor 2\ell/3\rfloor ! 9^{\ell/3}} \left(4K^2 \left(\ln \log_q g + O(1)\right)\right)^{\ell} \\
&\ll |\mathcal{C}_g| \Big( \frac{A}{KV} \Big)^{2\ell} \left( \frac{\ell}{e} \right)^{4\ell/3} \left( \frac{5}{3} \right)^{2\ell/3} (\ln \log_q g)^{\ell} \\
&\ll |\mathcal{C}_g| \exp\left( -\frac{V}{20A} \log V \right).
\end{aligned}$$

Now, if $\chi$ primitive such that $\mathcal{T}_1(\chi) \geq V_1$, then again applying Lemma 3.2 and Stirling's formula, we get

$$\begin{aligned}
\sum_{\chi\in\mathcal{C}_g} \left|\mathcal{T}_1(\chi)\right|^{2\ell} &\ll |\mathcal{C}_g| \frac{(\ell!)^2 5^{2\ell/3}}{\lfloor 2\ell/3\rfloor ! 9^{\ell/3}} \left( \sum_{d(P)\leq N_0} \frac{|h_1(P)|^2 |a(P)|^2}{|P|} \right)^{\ell} \\
&\ll |\mathcal{C}_g| \left( \frac{\ell \boldsymbol{\sigma}}{e} \right)^{\ell}
\end{aligned}$$

for any $\ell$ such that $2\ell N_0 \leq n$, which implies that $\ell \leq \frac{V}{A} \log_q g$. Markov's inequality gives us

$$\Psi_g(\boldsymbol{\theta}, V_1) \ll \frac{1}{V_1^{2\ell}} \left( \sum_{\chi\in\mathcal{C}_g} |\mathcal{T}_1(\chi)|^{2\ell} \right) \ll |\mathcal{C}_g| \left( \frac{\ell \boldsymbol{\sigma}}{e{V_1}^2} \right)^{\ell}.$$

It is now convenient to consider the case when $V \leq \frac{\boldsymbol{\sigma}^2}{K^3}$ and the case $V > \frac{\boldsymbol{\sigma}^2}{K^3}$ separately.

*Case* 1. Assume that $V \leq \frac{\boldsymbol{\sigma}^2}{K^3}$. We choose $\ell = \lfloor \frac{{V_1}^2}{\boldsymbol{\sigma}} \rfloor$. The definition of $A$ and this choice of $\ell$ implies that $\ell \leq \frac{V}{A} \log_q g$. In this case, we find that

$$\Psi_g(\boldsymbol{\theta}, V_1) \ll |\mathcal{C}_g| \ \exp\left( \ell \log \left( \frac{\ell \boldsymbol{\sigma}}{e{V_1}^2} \right) \right) \ll |\mathcal{C}_g| \ \exp\left( -\frac{{V_1}^2}{\boldsymbol{\sigma}} \right).$$

*Case* 2. Assume that $V > \frac{\boldsymbol{\sigma}^2}{K^3}$. We choose $\ell = \lfloor 10V \rfloor$. Again from the definition of $A$, it is easy to see that this choice $\ell$ satisfies $\ell \leq \frac{V}{A}\log_q g$. Notice that $V > \frac{\boldsymbol{\sigma}^2}{K^3}$, implies $\log V > 2\log\boldsymbol{\sigma} - 3\log K$. So, we have

$$A = 7K \quad \text{and} \quad {V_1}^2 = \frac{1}{49}V^2.$$

Hence, we conclude that

$$\Psi_g(\boldsymbol{\theta}, V_1) \ll |\mathcal{C}_g| \exp\left(10V\log\left(\frac{10V\boldsymbol{\sigma}}{e{V_1}^2}\right)\right).$$

Since

$$V > \frac{\boldsymbol{\sigma}^2}{K^3},$$

it follows that

$$\ln\boldsymbol{\sigma} < \frac{1}{2}\ln V + \frac{3}{2}\ln K.$$

Substituting this into the previous expression, we obtain

$$10V\ln\left(\frac{10V\boldsymbol{\sigma}}{eV_1^2}\right) \leq 10V\left(\ln(490) - 1 + \frac{3}{2}\ln K - \frac{1}{2}\ln V\right) = -5V\ln V + O_K(V).$$

Since $K$ is fixed, for sufficiently large $V$,

$$-5V\ln V + O_K(V) \leq -4V\ln V,$$

that gives

$$\Psi_g(\boldsymbol{\theta}, V_1) \ll |\mathcal{C}_g| \exp(-4V\ln V).$$

Therefore combining the above estimates, we deduce that

$$\Psi_g(\boldsymbol{\theta}, V) \ll |\mathcal{C}_g| \left\{\exp\left(-\tfrac{V}{20A}\ln V\right) + \exp\left(-\tfrac{{V_1}^2}{\boldsymbol{\sigma}}\right) + \exp\left(-4V\ln V\right)\right\}. \tag{3.6}$$

We now compute $V_1$ in each range of $V$ determined by the definition of $A$. If $\sqrt{\ln g} \leq V \leq \boldsymbol{\sigma}$, then

$$A = \frac{1}{2}\ln\boldsymbol{\sigma} \quad \text{and} \quad V_1 = V\Big(1 - \frac{12K}{\ln\boldsymbol{\sigma}}\Big).$$

So, for sufficiently large $g$, (3.6) implies that

$$\begin{aligned}\Psi_g(\boldsymbol{\theta}, V) &\ll |\mathcal{C}_g| \exp\left(-\tfrac{V^2}{\boldsymbol{\sigma}}\left(1 - \tfrac{12K}{\ln\boldsymbol{\sigma}}\right)^2\right) \\ &\ll |\mathcal{C}_g| \exp\left(-\tfrac{V^2}{\boldsymbol{\sigma}}\left(1 - \tfrac{24K}{\ln\boldsymbol{\sigma}}\right)\right).\end{aligned}$$

If $\boldsymbol{\sigma} \leq V \leq \frac{1}{25K}\boldsymbol{\sigma}\ln\boldsymbol{\sigma}$, then

$$A = \frac{\boldsymbol{\sigma}\log\boldsymbol{\sigma}}{2V} \quad \text{and} \quad V_1 = V\Big(1 - \tfrac{12KV}{\boldsymbol{\sigma}\ln\boldsymbol{\sigma}}\Big).$$

For this range of $V$, $\frac{\ln V}{\boldsymbol{\sigma}\ln\boldsymbol{\sigma}} > \frac{1}{\boldsymbol{\sigma}}$ and hence from (3.6) we obtain

$$\begin{aligned}\Psi_g(\boldsymbol{\theta}, V) &\ll |\mathcal{C}_g|\Bigg\{ \exp\left(-\frac{V^2\ln V}{\boldsymbol{\sigma}\ln\boldsymbol{\sigma}}\right) + \exp\left(-4V\ln V\right) \\ &\qquad + \exp\left(-\tfrac{V^2}{\boldsymbol{\sigma}}\left(1-\frac{12KV}{\boldsymbol{\sigma}\ln\boldsymbol{\sigma}}\right)^2\right)\Bigg\} \\ &\ll |\mathcal{C}_g|\exp\left(-\frac{V^2}{\boldsymbol{\sigma}}\left(1-\frac{24KV}{\boldsymbol{\sigma}\ln\boldsymbol{\sigma}}\right)\right).\end{aligned}$$

Finally, if $V > \frac{1}{25K}\boldsymbol{\sigma}\log\boldsymbol{\sigma}$, then

$$A = 7K \quad \text{and} \quad V_1 = \frac{V}{7}.$$

So from (3.6), we get that

$$\Psi_g(\boldsymbol{\theta}, V) \ll |\mathcal{C}_g|\exp\left(-\tfrac{V}{98K}\ln V\right).$$

Since for large $g$, $\boldsymbol{\sigma} = O_K(\ln g)$. In the range $\sqrt{\ln g} \le V \le 2026\boldsymbol{\sigma}$, we have

$$\exp\left(-\tfrac{V^2}{\boldsymbol{\sigma}}\left(1-\tfrac{24K}{\ln\boldsymbol{\sigma}}\right)\right) = \exp\left(-\tfrac{V^2}{\boldsymbol{\sigma}}\right)\exp\left(O_K\left(\frac{\ln g}{\ln\ln g}\right)\right) = \exp\left(-\tfrac{V^2}{\boldsymbol{\sigma}}\right)g^{\varepsilon}.$$

Combining these estimates with (3.6) for different range of $V$, we conclude that

$$\Psi_g(\boldsymbol{\theta}, V) \ll \begin{cases} |\mathcal{C}_g|g^{\varepsilon}\exp\left(-\tfrac{V^2}{\boldsymbol{\sigma}}\right) & \text{if } 3 \le V \le 2026\boldsymbol{\sigma}, \\ |\mathcal{C}_g|g^{\varepsilon}\exp\left(-4V\right) & \text{if } V > 2026\boldsymbol{\sigma}. \end{cases} \tag{3.7}$$

Substituting (3.7) into (3.2) completes the proof of Proposition 3.1.

## 4. Proof of Theorem 1.1

We consider the case $q \equiv 2 \pmod 3$. Using Lemmas 2.1 and 2.3, one can derive an estimate analogous to that in Theorem 1.1 for the Kummer extension.

4.0.1. *Set-up.* We closely follow the arguments from [20]. To find an upper bound for $S_g(\boldsymbol{\theta}^{(m)}, \boldsymbol{k}^{(m)})$, we can express it as

$$\sum_{\chi\in\mathcal{C}_g}\exp\left(2\sum_{j=1}^{m}k_j\ln\left|\mathcal{L}\left(\frac{e(\theta_j)}{q^{\frac{1}{2}+\alpha_j}},\chi\right)\right|\right).$$

Using Lemma 2.4, similar to eq. (3.3), we can provide an upper bound for the above quantity inside the exponent by a linear combination of Dirichlet polynomials. This means that, for $N \leq g+2$, we have

$$\ln\left(\left|\mathcal{L}\Big(\frac{e(\theta_1)}{q^{\frac{1}{2}+\alpha_1}},\chi\Big)\right|^{2k_1}\cdots\left|\mathcal{L}\Big(\frac{e(\theta_m)}{q^{\frac{1}{2}+\alpha_m}},\chi\Big)\right|^{2k_m}\right) = 2\,\mathrm{Re}\sum_{d(P)\leq N}\frac{\chi(P)h(P)a(P)}{|P|^{\frac{1}{2}+\frac{1}{N\ln q}}}$$

$$+\,2\,\mathrm{Re}\sum_{d(P)\leq N/2}\frac{\chi(P^2)h(P^2)a(P^2)}{2|P|^{1+\frac{2}{N\ln q}}}+\frac{2K(g+2)}{N}+O(K), \tag{4.1}$$

where for $P\in\mathcal{P}$, we define

$$h(P)=\sum_{j=1}^{m}k_j|P|^{-\alpha_j+\frac{2\pi i\theta_j}{\ln q}}, \tag{4.2}$$

and $a(P)=1-\frac{\deg P}{N}$. Certainly $a(P)\ll 1$, since $\deg(P)\geq 1$.

4.0.2. *Strategy to control Dirichlet polynomial over the family.* Instead of working with a single long Dirichlet polynomial, we decompose the sum over degrees of irreducible polynomials into $J+1$ dyadic intervals. This yields a collection of shorter Dirichlet polynomials supported on disjoint degree ranges.

Heuristically, the Dirichlet polynomials arising from distinct dyadic intervals are approximately independent. Consequently, their joint structure can be modeled by a *Brunching process*, analogous to the construction in [2], which is designed to capture the large values.

In contrast, when one works with a single long interval, this multiscale structure is not visible, and the underlying log-correlated behavior is not properly reflected.

Define $I_0 := (0,(g+2)t_0]$ and for $1\leq j\leq J$,

$$I_j=((g+2)t_{j-1},(g+2)t_j]$$

such that

$$t_j=\frac{e^{j-1}}{(\ln g)^2},\qquad 0\leq j\leq J.$$

Here, we choose $J=\left[\ln(\ln g)^2+\ln(t_J)\right]$ in such a way that $t_J=e^{-20K/(1-b)}$ is a small constant.

For each $1\leq i\leq j\leq J$, we denote

$$\mathcal{D}_{(i,j)}(\chi) := \mathcal{D}_{(i,j)}(\chi,g)=\sum_{\deg(P)\in I_i}\frac{a(P;j)\chi(P)h(P)}{|P|^{\frac{1}{2}+\frac{1}{(g+2)t_j\ln q}}} \tag{4.3}$$

where $a(P;j) = 1 - \frac{\deg P}{(g+2)t_j}$. For each $\chi \in \mathcal{C}_g$, the collection

$$\{\operatorname{Re}\mathcal{D}_{(i,j)}(\chi)\}_{1\le i\le j\le J}$$

may be viewed as forming an upper triangular matrix indexed by $(i,j)$. We now partition $\mathcal{C}_g$ by imposing size conditions on the entries in the last column of this matrix. More precisely, we define a subset

$$\mathcal{C}_g(J) = \left\{\chi \in \mathcal{C}_g : \left|\operatorname{Re}\mathcal{D}_{(i,J)}(\chi)\right| \le \gamma_i \ \forall\, 1 \le i \le J\right\},$$

where $\gamma_i = K[t_i^{-b}]$ for some $1/2 < b < 1$ that will be chosen later.

For $0 \le j \le J-1$, we define all other remaining possibilities as

$$\begin{aligned}\mathcal{C}_g(j) := \Big\{\chi \in \mathcal{C}_g : &\left|\operatorname{Re}\mathcal{D}_{(i,\ell)}(\chi)\right| \le \gamma_i \ \forall\, 1 \le i \le j, \quad \forall\, i \le \ell \le J \\ &\text{but } \left|\operatorname{Re}\mathcal{D}_{(j+1,\ell)}(\chi)\right| > \gamma_{j+1} \text{ for some } \ j+1 \le \ell \le J\Big\}.\end{aligned} \tag{4.4}$$

In this way, we can write

$$\mathcal{C}_g = \bigcup_{j=0}^{J} \mathcal{C}_g(j). \tag{4.5}$$

4.0.3. *Bounding exponent of Dirichlet polynomials over segments of $\mathcal{C}_g$.* We now state several auxiliary lemmas that are essential for completing the proof of the theorem. The first lemma provides a bound for the exponential of the sum of $\mathcal{D}_{(i,J)}(\chi)$'s appearing on the right-hand side of (4.1), when restricted to the subfamily $\mathcal{C}_g(J)$.

**Lemma 4.1.** *Assume that $g > \exp(4K^2)$. With the notations from (4.3) and (4.5), we have*

$$\sum_{\chi\in\mathcal{C}_g(J)} \exp\left(2\sum_{i=1}^{J}\operatorname{Re}\mathcal{D}_{(i,J)}(\chi)\right) \ll_{K,q} (\#\mathcal{C}_g)\cdot g^{\sum_{j=1}^m k_j^2} \prod_{1\le i<j\le m}\left(\min\left\{g, \frac{1}{|\theta_i-\theta_j|}\right\}\right)^{2k_ik_j}.$$

*Proof.* Using [23, Lemma 5.2], for any $\chi \in \mathcal{C}(J)$, we have

$$\exp\left(2\sum_{i=1}^{J}\operatorname{Re}\mathcal{D}_{(i,J)}(\chi)\right) \ll \prod_{i=1}^{J}\left(\sum_{0\le m\le[e^2\gamma_i]}\frac{(\operatorname{Re}\mathcal{D}_{(i,J)}(\chi))^m}{m!}\right)^2.$$

Therefore

$$\begin{aligned}\sum_{\chi\in\mathcal{C}_g(J)} \exp\left(\sum_{1\le j\le J} 2\operatorname{Re}\mathcal{D}_{(i,J)}(\chi)\right) &\ll \sum_{\chi\in\mathcal{C}_g(J)}\prod_{i=1}^{J}\left(\sum_{0\le m\le[e^2\gamma_i]}\frac{(\operatorname{Re}\mathcal{D}_{(i,J)}(\chi))^m}{m!}\right)^2 \\ &\le \sum_{\chi\in\mathcal{C}_g}\prod_{i=1}^{J}\left(\sum_{0\le m\le[e^2\gamma_i]}\frac{(\operatorname{Re}\mathcal{D}_{(i,J)}(\chi))^m}{m!}\right)^2.\end{aligned} \tag{4.6}$$

Note that, in the last line we sum over all $\chi \in \mathcal{C}_g$, rather than restricting to $\chi \in \mathcal{C}_g(J)$, using the positivity of the square of the inner sum.

Expanding first the $m$-th powers and then square, and moving forward the family of characters, the right hand side of (4.6) is equal to

$$\begin{aligned}
&\sum_{\widetilde{\alpha},\widetilde{\beta}} \left( \prod_{1\le i\le J} \frac{1}{\alpha_i!\beta_i!} \right) \sum_{\widetilde{P},\widetilde{Q}} C(\widetilde{P},\widetilde{Q}) \\
&\times \sum_{\chi\in\mathcal{C}_g} \prod_{1\le i\le J} \prod_{\substack{1\le r\le\alpha_i\\1\le s\le\beta_i}} \operatorname{Re}\left(\chi(P_{i,r})h(P_{i,r})\right) \operatorname{Re}\left(\chi(Q_{i,s})h(Q_{i,s})\right),
\end{aligned} \tag{4.7}$$

where $0 \le \alpha_i, \beta_i \le [e^2\gamma_i]$ with

$$\widetilde{\alpha} = (\alpha_1, \ldots, \alpha_J), \qquad \widetilde{\beta} = (\beta_1, \ldots, \beta_J).$$

Also here

$$\widetilde{P} = (P_{1,1}, \ldots, P_{1,\alpha_1}, \ldots, P_{J,1}, \ldots, P_{J,\alpha_J})$$

and

$$\widetilde{Q} = (Q_{1,1}, \ldots, Q_{1,\beta_1}, \ldots, Q_{J,1}, \ldots, Q_{J,\beta_J})$$

such that for $1 \le i \le J$,

$$P_{i,1}, \ldots, P_{i,\alpha_J},\, Q_{i,1}, \ldots, Q_{i,\beta_J} \in ((g+2)t_{i-1}, (g+2)t_i].$$

Lastly,

$$C(\widetilde{P},\widetilde{Q}) = \prod_{1\le i\le J} \prod_{\substack{1\le r\le\alpha_i\\1\le s\le\beta_i}} \frac{a(P_{i,r};J)}{|P_{i,r}|^{1/2+1/((g+2)t_J\ln q)}} \frac{a(Q_{i,s};J)}{|Q_{i,s}|^{1/2+1/((g+2)t_J\ln q)}}.$$

We use the identity $\operatorname{Re} z = \frac{z+\bar{z}}{2}$ to write the innermost sum of (4.7) as

$$\sum_{\chi\in\mathcal{C}_g} \prod_{1\le i\le J} \prod_{\substack{1\le r\le\alpha_i\\1\le s\le\beta_i}} \left( \frac{\chi(P_{i,r})h(P_{i,r}) + \overline{\chi(P_{i,r})h(P_{i,r})}}{2} \right) \left( \frac{\chi(Q_{i,s})h(Q_{i,s}) + \overline{\chi(Q_{i,s})h(Q_{i,s})}}{2} \right).$$

Using Lemmas 2.1 and 2.2, (4.7) is equal to

$$\begin{aligned}
&(\#\mathcal{C}_g) \prod_{1\le j\le J} \frac{1}{2^{\alpha_j+\beta_j}} \sum_{(\epsilon_{i,1},\ldots,\epsilon_{i,\alpha_i},\epsilon'_{i,1},\ldots,\epsilon'_{i,\beta_i})} \prod_{\substack{1\le r\le\alpha_i\\1\le s\le\beta_i}} h^{(\epsilon_{j,r})}(P_{j,r}) h^{(\epsilon'_{j,s})}(Q_{j,s}) \phi\Bigg( \prod_{1\le j\le J} \prod_{\substack{1\le r\le\alpha_j\\1\le s\le\beta_j}} P_{j,r}^{\epsilon_{j,r}} Q_{j,s}^{\epsilon'_{j,s}} \Bigg) \\
&\qquad + O\Bigg( q^{(1/2+\epsilon)g} \prod_{1\le j\le J} \left(\frac{K^2}{2}\right)^{\alpha_j+\beta_j} \prod_{1\le j\le J} \prod_{\substack{1\le r\le\alpha_j\\1\le s\le\beta_j}} |P_{j,r}Q_{j,s}|^{\epsilon} \Bigg),
\end{aligned} \tag{4.8}$$

where $(\epsilon_{i,1}, \ldots, \epsilon_{i,\alpha_i}, \epsilon'_{i,1}, \ldots, \epsilon'_{i,\beta_i})$, with each component equal to either $-1$ or $+1$. We adopt the convention that $h^{(1)}(P) = h(P)$ and $h^{(-1)}(P) = \overline{h(P)}$. Moreover, by (4.2) we have $|h(P)| \leq K$. Also here for any $f \in \mathbb{F}_q[t]$,

$$\phi(f) = \begin{cases} \prod_{P|f;2|d(P)} \frac{|P|}{2+|P|} & \text{if } f = \square, \\ 0 & \text{otherwise.} \end{cases}$$

We first evaluate the error term comes from the above expression. The contribution from the $O$-term in (4.8) is bounded above by

$$\begin{aligned} &\ll q^{(1/2+\epsilon)g} \prod_{1\leq i\leq J} \left( \sum_{0\leq\alpha_i\leq[e^2\gamma_i]} \frac{K^{\alpha_i}}{\alpha_i!} \left( \sum_{\deg(P)\in I_i} \frac{1}{|P|^{1/2-\epsilon}} \right)^{\alpha_i} \right)^2 \\ &\ll q^{(1/2+\epsilon)g} \prod_{1\leq i\leq J} q^{2e^2(1-2\epsilon)(g+2)t_i\gamma_i} \left( \sum_{0\leq\alpha_i\leq[e^2\gamma_i]} \frac{K^{\alpha_i}}{\alpha_i!} \right)^2 \\ &\ll q^{(1/2+\epsilon)g+4e^2(1-2\epsilon)(g+2)t_J\gamma_J} \exp(2KJ) \ll q^{(2/3+\epsilon)g} (\ln g)^{4K}, \end{aligned} \tag{4.9}$$

since we have chosen $J = \left[\ln(\ln g)^2 + \ln(t_J)\right]$ so that $t_J = e^{-20K/(1-b)}$.

Next, we turn to the main term in (4.8) and it is bounded above by

(4.10)

$$\begin{aligned} &\ll (\#\mathcal{C}_g) \times \prod_{1\leq i\leq J} \sum_{0\leq\alpha_i,\beta_i\leq[e^2\gamma_i]} \frac{1}{2^{\alpha_i+\beta_i}\alpha_i!\beta_i!} \sum_{(\epsilon_{i,1},\ldots,\epsilon_{i,\alpha_i},\epsilon'_{i,1},\ldots,\epsilon'_{i,\beta_i})} \\ &\times \sum_{\widetilde{P},\widetilde{Q}} \left| C(\widetilde{P},\widetilde{Q}) \prod_{\substack{1\leq r\leq\alpha_i \\ 1\leq s\leq\beta_i}} h^{(\epsilon_{j,r})}(P_{j,r}) h^{(\epsilon'_{j,s})}(Q_{j,s}) \right| \phi\Bigg( \prod_{\substack{1\leq r\leq\alpha_j \\ 1\leq s\leq\beta_j}} P_{j,r}^{\epsilon_{j,r}} Q_{j,s}^{\epsilon'_{j,s}} \Bigg) \\ &\ll (\#\mathcal{C}_g) \times \prod_{1\leq i\leq J} \sum_{0\leq\alpha_i,\beta_i\leq[e^2\gamma_i]} \frac{1}{2^{\alpha_i+\beta_i}\alpha_i!\beta_i!} \sum_{(\epsilon_{i,1},\ldots,\epsilon_{i,\alpha_i},\epsilon'_{i,1},\ldots,\epsilon'_{i,\beta_i})} \\ &\times \sum_{\substack{P_{i,1},\ldots,P_{i,\alpha_i},Q_{i,1},\ldots,Q_{i,\beta_i}\in I_i \\ \prod_{\substack{1\leq r\leq\alpha_i \\ 1\leq s\leq\beta_i}} P_{i,r}^{\epsilon_{i,r}} Q_{i,s}^{\epsilon'_{i,s}} = \square}} \prod_{\substack{1\leq r\leq\alpha_i \\ 1\leq s\leq\beta_i}} \frac{|a(P_{i,r};J)h^{(\epsilon_{i,r})}(P_{i,r})|}{|P_{i,r}|^{1/2+1/((g+2)t_J\log q)}} \frac{|a(Q_{i,s};J)h^{(\epsilon'_{j,s})}(Q_{j,s})|}{|Q_{i,s}|^{1/2+1/((g+2)t_J\log q)}}, \end{aligned}$$

where we have used that $\phi(f) \leq 1$ while $f = \square$.

For a fixed $1 \leq i \leq J$ and tuples $(\epsilon_{i,1}, \ldots, \epsilon_{i,\alpha_i}, \epsilon'_{i,1}, \ldots, \epsilon'_{i,\beta_i})$, we estimate the number of ways the following can occur:

$$P_{i,1}^{\epsilon_{i,1}} \cdots P_{i,\alpha_i}^{\epsilon_{i,\alpha_i}} Q_{i,1}^{\epsilon'_{i,1}} \cdots Q_{i,\beta_i}^{\epsilon'_{i,\beta_i}} = \square. \tag{4.11}$$

Note that if three irreducible polynomials of the form $(P, P, P)$ produce a cube, then the corresponding contribution is a Dirichlet polynomial of the form

$$E_i := \sum_{\deg(P)\in I_i} \frac{|a(P;J)|^3|h(P)|^3}{|P|^{3/2}} \ll K^3 \tag{4.12}$$

for some $1 \leq i \leq J$.

On the other hand, consider pairs of the form $(P, P^2)$. The contribution from such pairs is bounded by

$$M_i := \sum_{\deg(P)\in I_i} \frac{|a(P;J)|^2|h(P)|^2}{|P|}. \tag{4.13}$$

Thus the main contribution effectively comes from pairs of the form $(P, P^2)$. Furthermore, the contribution from $3k$-tuples with $k \in \mathbb{N}$ is also negligible and can absorb an occurrence of the quantity (4.13). More precisely, the contribution coming from $m \geq 2$ many pairs of the form $(P, P^2)$ is bounded above by $(M_i)^m$. If instead a triplet $(P, P, P)$ occurs, then the contribution is reduced to $(M_i)^{m-1}K^3$. This is negligible compared with $M_i$ provided that $\frac{K^3}{M} \ll 1$. Since, by (4.15) and (4.16), $M = \sum_i M_i \asymp K^2 \ln g$, the above condition follows from $g > \exp(4K^2)$.

Consequently, the main contribution comes from the number of ways

$$P_1 \cdots P_m\, Q_1^2 \cdots Q_m^2 = ⧈$$

such that $P_1, \ldots, P_m$ are distinct and $Q_1, \ldots, Q_m$ are so. Here $P_j$'s can pair ups with $Q_i$'s for $1 \leq i, j \leq m$. The number of such ways is exactly $m!$. In this case we may take $\alpha_i = \beta_i = m$.

Using the binomial expansion and the discussion above, we can bound the innermost sum on the right-hand side of the inequality (4.10) by

$$\ll \sum_{r\leq m} (M_i + E_i)^r \ll m m! M_i^m,$$

since $g > \exp(4K^2)$, where $M_i$ and $E_i$ are as in (4.13) and (4.12) respectively.

Based on the above analysis, we conclude that (4.7) is bounded by

$$\begin{aligned} (4.7) &\ll q^g \prod_{1\leq i\leq J} \sum_{m\leq[e^2\gamma_i]} \frac{m!\, m}{(m!)^2 4^m} \left( \sum_{d(P)\in I_i} \frac{|h(P)a(P;J)|^2}{|P|} \right)^m \\ &\ll q^g \exp\left( \sum_{\deg P\leq(g+2)t_J} \frac{|h(P)|^2}{|P|} \right) \\ &\ll q^g \cdot g^{k_1^2+k_2^2+\cdots+k_m^2} \prod_{1\leq i<j\leq m} \left( \min\left\{ g, \frac{1}{|\theta_i-\theta_j|} \right\} \right)^{2k_ik_j}, \end{aligned} \tag{4.14}$$

since $|a(P;J)| \ll 1$,

$$|h(P)|^2 = \sum_{j=1}^{m} k_j^2 |P|^{-2\alpha_j} + 2\sum_{i<j} k_i k_j |P|^{-\alpha_i-\alpha_j} \cos((\theta_i-\theta_j)\deg(P)), \tag{4.15}$$

and we know that [17, Appendix 9.1]

$$\sum_{n\le g} \frac{\cos(n\theta)}{n} = \ln\left(\min\left\{g, \tfrac{1}{|\theta|}\right\}\right) + O(1). \tag{4.16}$$

This completes the proof of the lemma. □

**Lemma 4.2.** *Consider the notations from* (4.3) *and* (4.5)*. Assume that* $g > \exp(4K^2)$*. Then we have*

$$meas(\mathcal{C}_g(0)) \ll (\#\mathcal{C}_g) \exp\left(-(\ln g)^2/4\right),$$

*and for any* $1 \le j \le J-1$*,*

$$\sum_{\chi\in\mathcal{C}_g(j)} \exp\left(2\sum_{i=1}^{j} \operatorname{Re}\mathcal{D}_{(i,j)}(\chi)\right)$$

$$\ll (\#\mathcal{C}_g)\exp\left(-t_{j+1}^{-1}\ln(1/t_{j+1})/13\right)\cdot g^{k_1^2+k_2^2+\cdots+k_m^2} \prod_{1\le i<j\le m}\left(\min\left\{g, \frac{1}{|\theta_i-\theta_j|}\right\}\right)^{2k_ik_j}.$$

*Proof.* From the definition of $\mathcal{C}_g(j)$ in (4.4), we can apply Rankin's trick to obtain that for any $M>0$,

$$\sum_{\chi\in\mathcal{C}_g(j)} \exp\left(2\sum_{i=1}^{j} \operatorname{Re}\mathcal{D}_{(i,j)}(\chi)\right) \le \sum_{\ell=j+1}^{J} \sum_{\substack{\chi\in\mathcal{C}_g \\ |\operatorname{Re}\mathcal{D}_{(i,j)}(\chi)|\le\gamma_i \forall 1\le i\le j \\ \text{but } |\operatorname{Re}\mathcal{D}_{(j+1,\ell)}(\chi)|>\gamma_{j+1}}} \prod_{1\le i\le j} \exp\left(2\operatorname{Re}\mathcal{D}_{(i,j)}(\chi)\right)$$

(4.17)

$$\ll (\gamma_{j+1})^{-2M} \sum_{\chi\in\mathcal{C}_g} \prod_{i=1}^{J}\left(\sum_{0\le m\le[e^2\gamma_i]} \frac{(\operatorname{Re}\mathcal{D}_{(i,J)}(\chi))^m}{m!}\right)^2 \left(\operatorname{Re}\mathcal{D}_{(j+1,\ell)}(\chi)\right)^{2M},$$

where we applied [23, Lemma 5.2] again similar way as of (4.6).

Since the length of the Dirichlet polynomial must be less than $(g+2)$, we are forced to impose the restriction

$$2(g+2)t_{j+1}M \le g+2 \implies M \le \left[\frac{1}{t_{j+1}}\right].$$

Accordingly, we choose $M = \left[\frac{1}{6t_{j+1}}\right]$. The sum on the right-hand side of the above inequality can be estimated in exactly the same way as in the proof of

Lemma 4.1, mainly by following the estimates of (4.7). Therefore, the inner sum of (4.17) is bounded above by

$$\begin{aligned}&\ll (\#\mathcal{C}_g) \prod_{1\le i\le j} \sum_{m\le[e^2\gamma_i]} \frac{m!\,m}{(m!)^2 4^m}\Bigg(\sum_{\deg(P)\in I_i} \frac{|h(P)|^2}{|P|}\Bigg)^m \\ &\times M^3(M!)\Bigg(\sum_{\deg(P)\in I_{j+1}} \frac{|h(P)|^2}{|P|}\Bigg)^M + O\left(q^{(7/8+\epsilon)g}(\ln g)^{2K}\right) \\ &\ll (\#\mathcal{C}_g)\exp\left(\sum_{\deg(P)\le(g+2)t_j} \frac{|h(P)|^2}{|P|}\right)\left(\frac{1}{6t_{j+1}}\sum_{\deg(P)\in I_{j+1}} \frac{|h(P)|^2}{|P|}\right)^{\left[\frac{1}{6t_{j+1}}\right]} \\ &\qquad + O\left(q^{(7/8+\epsilon)g}(\ln g)^{2K}\right),\end{aligned}$$

where we separated the variables $m$ and $M$ using the inequality $\phi(fg) \le \phi(f)\phi(g)$ for any $f, g \in \mathbb{F}_q[t]$, together with the bound $\phi(f) \le 1$, and Stirling's approximation has also been applied. Note that we used the fact that the two Dirichlet polynomials supported on distinct set of irreducible polynomials. Furthermore, the error term is estimated in exactly the same way as in (4.9), with an additional contribution to the exponent arising from the choice of $M$. Although the error term can be further optimized, the present bound is sufficient for our purposes.

Finally, putting altogether in (4.17) and following the computations carried out in (4.14), we can conclude that

$$\begin{aligned}&\sum_{\chi\in\mathcal{C}_g(j)} \exp\left(2\sum_{i=1}^{j} \operatorname{Re}\mathcal{D}_{(i,j)}(\chi)\right) \ll (I-j)\,(\#\mathcal{C}_g)\cdot\, g^{k_1^2+k_2^2+\cdots+k_m^2} \\ &\times \prod_{1\le i<j\le m}\left(\min\left\{g, \frac{1}{|\theta_i-\theta_j|}\right\}\right)^{2k_ik_j}\left(\frac{t_{j+1}^{2b-1}}{6K^2}\sum_{\deg P\in I_j} \frac{|h(P)|^2}{|P|}\right)^{\left[\frac{1}{6t_{j+1}}\right]},\end{aligned} \tag{4.18}$$

provided that $2b-1>0$. Here we choose $b = 3/4$.

**Case 1:** Assume that $j = 0$. In this case, $t_0 = 1/(e(\ln g)^2)$, $I_0 = (0, (g+2)t_0]$, and $J \ll \ln_2 g$. Also by Prime polynomial theorem,

$$\sum_{\deg(P)\le(g+2)t_1} 1/|P| \le \ln g.$$

It follows that

$$\begin{aligned}\operatorname{meas}(\mathcal{C}_g(0)) &\ll (\#\mathcal{C}_g)\,(\ln_2 g)\exp(K^2\ln g)\exp\left(\frac{1}{6t_1}\ln\left((t_1^{1/2}\ln g)/6\right)\right)\\ &\ll (\#\mathcal{C}_g)\,(\ln_2 g)\exp\left(K^2\ln g-(\ln 6)(\ln g)^2/6\right)\\ &\ll (\#\mathcal{C}_g)\exp\left(-(\ln g)^2/4\right)\end{aligned}$$

provided that $g>\exp(4K^2)$.

**Case 2:** Let $1\le j\le I-1$. In this case, $I-j\ll \ln_2 g$, and again by prime polynomial theorem, we have

$$\sum_{\deg P\in I_j}\frac{1}{|P|}\le \sum_{(g+2)t_{j-1}<n\le (g+2)t_j} 1/n=\ln t_j-\ln t_{j-1}+o(1)=1+o(1)\le 2.$$

Therefore, we conclude from (4.18) that

$$\begin{aligned}&\sum_{\chi\in\mathcal{C}_g(j)}\exp\left(2\operatorname{Re}\mathcal{D}_{(j,j)}(\chi)\right)\ll (\#\mathcal{C}_g)\cdot\, g^{k_1^2+k_2^2+\cdots+k_m^2}\\ &\times\prod_{1\le i<j\le m}\left(\min\left\{g,\frac{1}{|\theta_i-\theta_j|}\right\}\right)^{2k_ik_j}\exp\left(-t_{j+1}^{-1}\ln(1/t_{j+1})/13\right).\end{aligned}$$

This completes the proof. $\square$

**Lemma 4.3.** *For $1\le j\le J$, the estimates in Lemmas 4.1 and 4.2 remain true if we replace the Dirichlet polynomial in the exponents by*

$$\mathcal{D}_{(j,j)}(\chi)+\sum_{\deg(P)\le (g+2)t_J/2}\frac{\chi(P^2)h(P^2)a(P^2;J)}{2|P|^{1+2/((g+2)t_J\ln q)}}.$$

*Proof.* We first show, by a modification of the proof of Lemma 4.1, that

$$\begin{aligned}&\sum_{\chi\in\mathcal{C}_g(J)}\exp\left(2\operatorname{Re}\left(\mathcal{D}_{(J,J)}(\chi)+\sum_{\deg(P)\le(g+2)t_J/2}\frac{\chi(P^2)h(P^2)a(P^2;J)}{2|P|^{1+2/((g+2)t_J\ln q)}}\right)\right)\\ &\ll_{K,q}(\#\mathcal{C}_g)\cdot\, g^{k_1^2+k_2^2+\cdots+k_m^2}\prod_{1\le i<j\le m}\left(\min\left\{g,\frac{1}{|\theta_i-\theta_j|}\right\}\right)^{2k_ik_j}.\end{aligned}\tag{4.19}$$

The modification of Lemma 4.2 to accommodate this new exponent can be carried out in a similar manner.

To prove (4.19), for $0\le m\le\frac{(g+2)t_J}{2}$, we introduce the set

$$\begin{aligned}\mathcal{B}_m(J):=\Big\{\chi\in\mathcal{C}_g\,:\,&|\operatorname{Re}\mathcal{P}_m(\chi)|>\alpha^{-m}\text{ but}\\ &|\operatorname{Re}\mathcal{P}_r(\chi)|\le\alpha^{-r}\text{ for all }m+1\le r\le\frac{(g+2)t_J}{2}\Big\},\end{aligned}$$

where $\alpha > 1$ are parameters to be chosen later, and decomposing the additional exponent as

$$\mathcal{P}_m(\chi) := \sum_{P \in \mathcal{P}_m} \frac{\chi(P^2)h(P^2)a(P^2;J)}{2|P|^{1+2/((g+2)t_J \ln q)}}.$$

Note that,

$$\mathcal{P}(\chi) := \sum_{m \le \frac{(g+2)t_J}{2}} \mathcal{P}_m(\chi) = \sum_{\deg(P) \le (g+2)t_J/2} \frac{\chi(P^2)h(P^2)a(P^2;J)}{2|P|^{1+2/((g+2)t_J \ln q)}}.$$

First we see that if $\chi \notin \mathcal{B}_m(J)$ for any $m \ge 1$, then $|\operatorname{Re}\mathcal{P}_r(\chi)| \le \alpha^{-r}$ for all $r$, and then summing over these $r$, we have $\mathcal{P}(\chi) \ll 1$, and so the part of $\sum_{\chi \in \mathcal{C}_g(J)}$ corresponding to such so called 'good' $\chi$ can be bounded exactly same way as in Lemma 4.1. To compute 'bad' $\chi$'s for which $\chi \in \mathcal{B}_m(J)$, we break the analysis into two parts depending on the size of $m$.

**Case 1:** First we consider that $m > [(4\log_q g)/5]$. In this case, we first show that while $q > \alpha^2$,

$$\operatorname{meas}(\mathcal{B}_m(J)) \ll q^g \exp\left(-m \ln(q/\alpha^2)\right). \tag{4.20}$$

From the construction of the set $\mathcal{B}_m(J)$, we have

$$\operatorname{meas}(\mathcal{B}_m(J)) \ll \alpha^{2m} \sum_{\chi \in \mathcal{C}_g} (\operatorname{Re}\mathcal{P}_m(\chi))^2.$$

Again we can repeat the arguments of the proof of Lemma 4.1, mainly (4.9) and (4.14) to obtain

$$\begin{aligned}\operatorname{meas}(\mathcal{B}_m(J)) &\ll q^g \alpha^{2m} \sum_{P \in \mathcal{P}_m} \frac{|h(P)|^2 \phi(P)}{|P|^2} + O\left(q^{(1/2+\epsilon)} \left(\sum_{P \in \mathcal{P}_m} \frac{|P|^\epsilon |h(P)|}{|P|}\right)^2\right) \\ &\ll q^g \alpha^{2m} K^2 q^{-m} + O\left(q^{(1/2+\epsilon)} \left(\frac{K\alpha^m q^{\epsilon m}}{m}\right)^2\right) \\ &\ll q^g \exp(-m\ln(q/\alpha^2)) + O\left(q^{(1/2+\epsilon)g + 2(1+\epsilon)m}\right) \\ &\ll q^g \exp(-m\ln(q/\alpha^2)) + O\left(q^{(2/3+\epsilon)g}\right),\end{aligned}$$

provided that $\alpha^2 < q$, and given that $m \le (g+2)t_J/2$. This completes the claim (4.20).

On the other hand, if $\chi \in \mathcal{B}_m(J)$, then by another application of the prime polynomial theorem, we obtain that

$$\operatorname{Re}\mathcal{P}(\chi) \le \sum_{\ell \le m} \frac{K}{2\ell} + \sum_{\ell=m+1}^{(g+2)t_J/2} \frac{K}{\alpha^\ell} \le \frac{K}{2}(\ln m + \gamma + 1) + \frac{K}{\alpha^{m+1}}(1 - 1/\alpha)^{-1}.$$

This lead us to obtain by using (4.20) that

$$\sum_{\chi\in\mathcal{B}_m(J)} \exp(4\,\mathrm{Re}\,\mathcal{P}(\chi)) \ll_K m^{2K}\mathrm{meas}(\mathcal{B}_m(J)) \ll_K q^g m^{2K}\exp(-m\ln(q/\alpha^2)),$$

where we have use the inequality that $\exp(\frac{4K}{\alpha^{m+1}}(1-1/\alpha)^{-1}) \leq \exp(20K)$ with the choice $\alpha = q^{1/10}$ and $q \geq 5$.

Combining the above estimates with Lemma 4.1, and applying the Cauchy–Schwarz inequality, we finally obtain

$$\begin{aligned}
&\sum_{\chi\in\mathcal{C}_g(J)\cap\mathcal{B}_m(J)} \exp(2\,\mathrm{Re}(\mathcal{D}_{(J,J)}(\chi)+\mathcal{P}(\chi)))\\
&\ll \Bigg(\sum_{\chi\in\mathcal{C}_g(J)} \exp(4\,\mathrm{Re}\,\mathcal{D}_{(J,J)}(\chi))\Bigg)^{1/2}\Bigg(\sum_{\chi\in\mathcal{B}_m(J)} \exp(4\,\mathrm{Re}\,\mathcal{P}(\chi))\Bigg)^{1/2}\\
&\ll_K q^g \cdot g^{k_1^2+k_2^2+\cdots+k_m^2} \prod_{1\leq i<j\leq m}\left(\min\left\{g, \frac{1}{|\theta_i-\theta_j|}\right\}\right)^{2k_ik_j} g^{K^2}m^{2K}\exp(-m\ln(q/\alpha^2)).
\end{aligned} \tag{4.21}$$

Summing over $m$ from $[(4\log_q g)/5]$ to $(g+2)t_J/2$ and applying Rankin's trick, the desired result follows. The factor $g^{K^2}$ is absorbed by the tail of the $m$-sum, together with the bound $\sum_{m\geq 1} m^K x^m \leq \frac{x\,K!}{(1-x)^{K+1}}$.

**Case 2:** Assume that $m \leq [(4\log_q g)/5]$. This means it will only capture $\chi \in \mathcal{C}_g(J)\cap\mathcal{B}_m(J)$ such that $0 \leq m \leq [(4\log_q g)/5]$. We first trivially bound that if $\chi\in\mathcal{B}_m(J)$,

$$\mathrm{Re}\Bigg(\sum_{\deg(P)\leq m+1} \frac{\chi(P)h(P)a(P;J)}{|P|^{1/2+1/((g+2)t_J\ln q)}} + \mathcal{P}(\chi)\Bigg) \ll_q Kq^{m/2}.$$

We introduce this truncation of $\mathcal{D}_{(J,J)}(\chi)$ in order to compute moments of two Dirichlet polynomials supported on disjoint sets of irreducibles, which will appear below.

Using this, together with Rankin's trick to handle the tail arising from the restriction to the set $\mathcal{B}_m(J)$, we obtain

$$\sum_{\chi \in \mathcal{C}_g(J)\cap \mathcal{B}_m(J)} \exp(2\,\mathrm{Re}(\mathcal{D}_{(J,J)}(\chi)+\mathcal{P}(\chi)))$$
$$\ll \exp\left(O(Kq^{m/2})\right) \sum_{\chi \in \mathcal{C}_g(J)\cap \mathcal{B}_m(J)} \exp\left(2\,\mathrm{Re} \sum_{m+1<\deg(P)\leq (g+2)t_J} \frac{\chi(P)h(P)a(P;J)}{|P|^{1/2+1/((g+2)t_J \ln q)}}\right)$$
$$\ll \exp\left(O(Kq^{m/2})\right) \alpha^{2M}$$
$$\times \sum_{\chi\in\mathcal{C}_g(J)} (\mathrm{Re}\,\mathcal{P}_m(\chi))^{2M} \exp\left(2\,\mathrm{Re} \sum_{m+1<\deg(P)\leq (g+2)t_J} \frac{\chi(P)h(P)a(P;J)}{|P|^{1/2+1/((g+2)t_J \ln q)}}\right),$$

where $M$ will be chosen appropriately at a later stage.

At this point, we are ready to apply the same strategy as in the proof of Lemma 4.1. Arguing similarly to the estimate in (4.17), and noting that the Dirichlet polynomials $\mathrm{Re}\,\mathcal{P}_m(\chi)$ and the one appearing in the exponent are supported on disjoint sets of irreducible, we deduce that the above quantity is bounded above by

$$\ll \exp\left(O(Kq^{m/2})\right) q^g \alpha^{2M} M^3 M! \exp\left(\sum_{m+1<d(P)\leq (g+2)t_J} \frac{|h(P)|^2}{|P|}\right)\left(\sum_{P\in\mathcal{P}_m} \frac{|h(P)|^2}{|P|^2}\right)^M$$
$$+ O\left(q^{(1/2+\epsilon)g+2e^2(g+2)t_J\gamma_J+2\epsilon m M} \exp\left(O(q^{m/2})\left(K^2\alpha^m\right)^{2M}\right)\right)$$
$$\ll \exp\left(O(Kq^{m/2})\right) q^g \exp\left(\sum_{d(P)\leq (g+2)t_J} \frac{|h(P)|^2}{|P|}\right) \exp\left(M \ln\left(K^2 M (\alpha^2/q)^m\right)\right)$$
$$+ O\left(q^{(2/3+\epsilon)g}\right)$$
$$\ll \exp\left(O(Kq^{m/2})\right) q^g \exp\left(\sum_{d(P)\leq (g+2)t_J} \frac{|h(P)|^2}{|P|}\right) \exp\left(-q^{3m/4}\right),$$

where we have picked up $M := q^{3m/4}$ with $m \leq [(4\log_q g)/5]$, $\alpha = q^{1/10}$ and the condition that $g > \exp(4K^2)$. Summing over the $m$ from 1 to $[(4\log_q g)/5]$ yields the desired bound as in (4.21), thereby completing the proof of the lemma. □

4.0.4. *Proof of Theorem 1.1.* We now ready to prove the theorem by using Lemmas from section 4.0.3. From the partition (4.5) of the the set $\mathcal{C}_g$, it is

enough to show that

$$\sum_{j=0}^{J-1} \sum_{\chi \in \mathcal{C}_g(j)} \exp\left(2\sum_{\ell=1}^m k_\ell \ln\left|\mathcal{L}\Big(\frac{e(\theta_\ell)}{q^{\frac{1}{2}+\alpha_\ell}}, \chi\Big)\right|\right) + \sum_{\chi \in \mathcal{C}_g(J)} \exp\left(2\sum_{\ell=1}^m k_\ell \ln\left|\mathcal{L}\Big(\frac{e(\theta_\ell)}{q^{\frac{1}{2}+\alpha_\ell}}, \chi\Big)\right|\right)$$
$$\ll_{K,q} (\#\mathcal{C}_g) \cdot g^{k_1^2+k_2^2+\cdots+k_m^2} \prod_{1 \le i < j \le m} \left(\min\left\{g, \frac{1}{|\theta_i - \theta_j|}\right\}\right)^{2k_ik_j}.$$

Using Lemma 4.3 (specially (4.19)), and eq. (4.1) with $N = (g+2)t_J$, we have

$$\sum_{\chi \in \mathcal{C}_g(J)} \exp\left(2\sum_{\ell=1}^m k_\ell \ln\left|\mathcal{L}\Big(\frac{e(\theta_\ell)}{q^{\frac{1}{2}+\alpha_\ell}}, \chi\Big)\right|\right)$$
$$\ll_K e^{2K/t_J} (\#\mathcal{C}_g) \cdot g^{k_1^2+k_2^2+\cdots+k_m^2} \prod_{1 \le i < j \le m} \left(\min\left\{g, \frac{1}{|\theta_i - \theta_j|}\right\}\right)^{2k_ik_j},$$

where we use the choice of $t_J = e^{-80K}$. Similarly, for $1 \le j \le J-1$, again Lemma 4.3 and eq. (4.1) with $N = (g+2)t_j$, we have

$$\sum_{\chi \in \mathcal{C}_g(j)} \exp\left(2\sum_{j=1}^m k_j \ln\left|\mathcal{L}\Big(\frac{e(\theta_j)}{q^{\frac{1}{2}+\alpha_j}}, \chi\Big)\right|\right) \ll (\#\mathcal{C}_g) \exp\big(2K/t_j - t_{j+1}^{-1} \ln(1/t_{j+1})/13\big)$$
$$\times g^{k_1^2+k_2^2+\cdots+k_m^2} \prod_{1 \le i < j \le m} \left(\min\left\{g, \frac{1}{|\theta_i - \theta_j|}\right\}\right)^{2k_ik_j}.$$

Since $t_{j+1} \le t_J = \exp(80K)$, and $t_{j+1} = et_j$,

$$\exp\big(2K/t_j - t_{j+1}^{-1} \ln(1/t_{j+1})/13\big) = \exp\left(\frac{2K}{t_j} - \frac{\ln(1/t_{j+1})}{13et_j}\right) \le \exp(-1/(4t_j)).$$

Summing over $1 \le j \le J-1$, we get the required bound. Finally, for $j = 0$, Using Lemma 4.2, and bound from Section 3, we have

$$\sum_{\chi \in \mathcal{C}_g(0)} \exp\left(2\sum_{j=1}^m k_j \ln\left|\mathcal{L}\Big(\frac{e(\theta_j)}{q^{\frac{1}{2}+\alpha_j}}, \chi\Big)\right|\right)$$
$$\ll (\mathrm{meas}(\mathcal{C}_g(0)))^{1/2} \left(\sum_{\chi \in \mathcal{C}_g} \exp\left(4\sum_{j=1}^m k_j \ln\left|\mathcal{L}\Big(\frac{e(\theta_j)}{q^{\frac{1}{2}+\alpha_j}}, \chi\Big)\right|\right)\right)^{1/2}$$
$$\ll q^g g^{k_1^2+k_2^2+\cdots+k_m^2} \prod_{1 \le i < j \le m} \left(\min\left\{g, \frac{1}{|\theta_i - \theta_j|}\right\}\right)^{2k_ik_j} g^{K^2+1/2} e^{-(\ln g)^2/2}$$
$$\ll_K q^g g^{k_1^2+k_2^2+\cdots+k_m^2} \prod_{1 \le i < j \le m} \left(\min\left\{g, \frac{1}{|\theta_i - \theta_j|}\right\}\right)^{2k_ik_j}$$

provided that $g > \exp(4K^2)$. This completes the proof of the theorem.

## 5. Proof of Corollary

1.2 Let $\Gamma_{1/g}$ be the circle in the complex plane whose center is at origin and radius is $\frac{1}{g}$. By Cauchy's integral formula, we have

$$\mathcal{L}^{(\ell)}\left(q^{-1/2}, \chi\right) = \frac{\ell!}{2\pi i} \oint_{\Gamma_{1/g}} L\left(\frac{1}{2}+s, \chi\right) \frac{ds}{s^{\ell+1}}.$$

Note that if $s = \alpha - i\frac{2\pi\theta}{\log q}$ then we can express

$$L\left(\frac{1}{2}+s, \chi\right) = \mathcal{L}\left(\frac{e(\theta)}{q^{\alpha+1/2}}, \chi\right).$$

Now by applying Hölder inequality, we obtain

$$\left|\oint_{\Gamma_{1/g}} L\left(\frac{1}{2}+s, \chi\right) \frac{ds}{s^{\ell+1}}\right| \le \left(\oint_{\Gamma_{1/g}} \left|L\left(1/2+s, \chi\right)\right|^{2k} |ds|\right)^{\frac{1}{2k}} \left(\oint_{\Gamma_{1/g}} |s|^{\frac{-2k(\ell+1)}{2k-1}} |ds|\right)^{\frac{2k-1}{2k}}.$$

Therefore,

$$\begin{aligned}\sum_{\chi\in\mathcal{C}_g} \left|\mathcal{L}^{(\ell)}\left(q^{-1/2}, \chi\right)\right|^{2k} &\ll \left(\frac{\ell!}{2\pi}\right)^{2k} \left(\frac{2\pi}{g}\right)^{2k-1} g^{2k(\ell+1)} \sum_{\chi\in\mathcal{C}_g} \oint_{\Gamma_{1/g}} \left|L\left(1/2+s, \chi\right)\right|^{2k} |ds| \\ &\ll \left(\frac{\ell!}{2\pi}\right)^{2k} \left(\frac{2\pi}{g}\right)^{2k} g^{2k(\ell+1)} \max_{|\theta|\le\frac{1}{g}} \sum_{\chi\in\mathcal{C}_g} \left|L\left(\frac{1}{2}+s, \chi\right)\right|^{2k}.\end{aligned}$$

As a direct application of Theorem 1.1, we get

$$\frac{1}{|\mathcal{C}_g|} \sum_{\chi\in\mathcal{C}_g} \left|\mathcal{L}\left(\frac{e(\theta)}{q^{\alpha+1/2}}, \chi\right)\right|^{2k} \ll g^{k^2}.$$

Using this upper bound to the above inequality, we conclude that

$$\frac{1}{|\mathcal{C}_g|} \sum_{\chi\in\mathcal{C}_g} \left|\mathcal{L}^{(\ell)}\left(\frac{e(\theta)}{q^{\alpha+1/2}}, \chi\right)\right|^{2k} \ll g^{k^2+2k\ell}.$$

## Acknowledgement

During the preparation of this work, P.D. was supported by the Australian Research Council Grant DP230100534. S.D. was supported by NBHM Fellowship, DAE, Govt. of India (Ref. no 0204/1 (19)/2022/RbD-11/1226) during the course of this work at Indian Statistical Institute, Kolkata. S.D. also gratefully acknowledges the hospitality and financial support of the Max Planck Institute for Mathematics, Bonn, during the preparation of this manuscript. Part of this work was carried out while G.M. was visiting the Centre de recherches mathématiques (CRM) under the thematic program: Universal Statistics in Number Theory. G.M. gratefully acknowledges the excellent working conditions and hospitality provided by the Centre de Recherches Mathématiques. Currently, G.M. is supported by the National Board for Higher Mathematics (NBHM) Post-doctoral grant 0204/21(24)/2025-R&D-II/16322.

School of Mathematics and Statistics, University of New South Wales, Sydney NSW 2052, Australia
*Email address*: darbarpranendu100@gmail.com

Independent Researcher, India.
*Email address*: sampa.math@gmail.com

Stat-Math Unit, Indian Statistical Institute, Kolkata 700108, India
*Email address*: g.gopaltamluk@gmail.com